\documentclass[aps,prd,amsmath,amssymb]{revtex4-1}
\usepackage[T1]{fontenc}
\usepackage[latin9]{inputenc}
\usepackage{subfig}
\usepackage{graphicx}
\usepackage{float}
\usepackage{mathrsfs}
\usepackage{amscd}
\usepackage{natbib}
 \usepackage[unicode=true,pdfusetitle,
  bookmarks=true,bookmarksnumbered=false,bookmarksopen=false,
  breaklinks=true,pdfborder={0 0 1},backref=false,colorlinks=true]
  {hyperref}

  \newcommand*\cv[2]{
      \begin{bmatrix}#1\\#2\end{bmatrix}
  }
  \newcommand*\rv[2]{
      \begin{bmatrix}#1 & #2\end{bmatrix}
  }

  \newcommand*\tbt[4]{
      \begin{bmatrix}#1 & #2 \\ #3 & #4\end{bmatrix}
  }

\begin{document}
\title{Honey from the Hives: A Theoretical and Computational Exploration of Combinatorial Hives}
\author{John Lombard}
\email[]{jml448@uw.edu}
\affiliation{University of Washington, \\ Seattle, Washington, 98195, USA}
\date{\today}

\begin{abstract}

In the first half of this manuscript, we begin with a brief review of combinatorial hives as introduced by Knutson and Tao, and focus on a conjecture by Danilov and Koshevoy for generating such a hive from Hermitian matrix pairs through an optimization scheme. We examine a proposal by Appleby and Whitehead in the spirit of this conjecture and analytically elucidate an obstruction in their construction for guaranteeing hive generation, while detailing stronger conditions under which we can produce hives with almost certain probability. We provide the first mapping of this prescription onto a practical algorithmic space that enables us to produce affirming computational results and open a new area of research into the analysis of the random geometries and curvatures of hive surfaces from select matrix ensembles.

The second part of this manuscript concerns Littlewood-Richardson coefficients and methods of estimating them from the hive construction. We illustrate experimental confirmation of two numerical algorithms that we provide as tools for the community: one as a rounded estimator on the continuous hive polytope volume following a proposal by Narayanan, and the other as a novel construction using a coordinate hit-and-run on the hive lattice itself. We compare the advantages of each, and include numerical results on their accuracies for some tested cases.

\end{abstract}

\maketitle

\section{What are Combinatorial Hives?}

In 1934, D. E. Littlewood and A. R. Richardson formulated a combinatorial rule for the multiplication of Schur polynomials \cite{Littlewood1934}. With the following decades producing a wide variety of proofs of this rule, a vast scope of applications throughout combinatorics and representation theory were discovered, with natural numbered Littlewood-Richardson coefficients (LRC) $c^\lambda_{\mu \nu}$ appearing in such examples as the decomposition of the tensor product of two Schur modules or the symmetries of Young tableaux \cite{Leeuwen2001}.

With respect to their use as the multiplicities in decompositions of tensor products on $GL_n(\mathbb{C})$ for triples of dominant weights $(\mu,\nu,\lambda)$, Knutson and Tao (KT) introduce the `hive' combinatorial model in their proof of the saturation conjecture which implies that satisfying a particular series of linear inequalities is actually a sufficient condition for the corresponding LRC to be positive \cite{Knutson1999}.

Specifically, let the spectra $\sigma$ of a symmetric non-negative definite matrix $M$ be given by $\sigma(M)$ such that the eigenvalues are listed in weakly decreasing order. We define the LR-Cone as the subspace of tuples $\sigma(M,N,L) \equiv (\mu,\nu,\lambda)$ such that $M+N = L$. The proof demonstrates that if $(\mu,\nu,\lambda)$ is an integer point in the LR-Cone, then the corresponding structure constant $c^\lambda_{\mu \nu}$ is greater than $0$.

\begin{figure}[H]
  \centerline{\includegraphics[scale = .55]{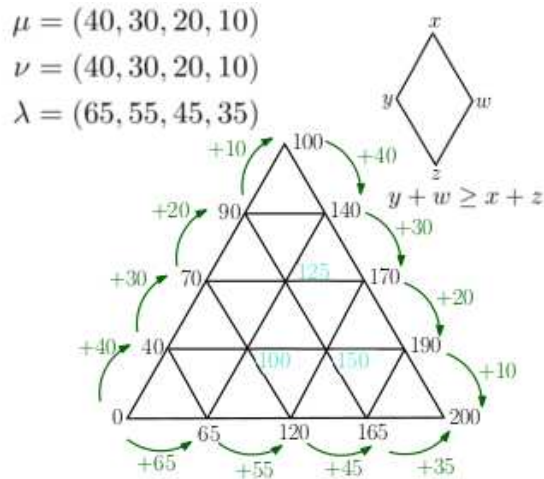}}
  \caption{An example hive boundary setup and a possible interior labeling. A graphical representation of a rhombus inequality is shown for a vertical rhombus. \cite{Narayanan2014}.}\label{fig:HiveExample}
\end{figure}

KT cast this condition into a geometric structure as follows: given a non-negative integer valued weakly decreasing 3-tuple of vectors $(\mu,\nu,\lambda)$ each of length $n$, one can use these vectors to construct the boundary values of a triangular hive, illustrated in Fig. \ref{fig:HiveExample}, that maps onto the previous condition provided that
\begin{equation}
  \|\mu\|_1 + \|\nu\|_1 = \|\lambda\|_1 \, . \label{eq:saturationCondition}
\end{equation}

The triangle boundaries are built out of the cumulative sums of $\mu$ on the left edge starting at $0$ on the left corner, $\mu + \nu$ continuing from the top into the right corner, and $\lambda$ working from the left corner to the right corner, with Eq. \ref{eq:saturationCondition} enforcing that the labels meet appropriately at the right corner. One can then tessellate the interior with equilateral triangles according to the lattice points of the boundary. If one assigns integer labels to the bulk interior points such that for every minimal rhombus in the triangulation
\begin{equation}
  \sum{\text{(Labels at the Obtuse Angles)}} \geq \sum{\text{(Labels at the Acute Angles)}} \,, \label{eq:rhombusAngles}
\end{equation}
then the resulting object is a member of the set $HIVE(\mu,\nu,\lambda)$. This condition enforces that the piecewise linear extension to the interior of the boundary triangle is a concave function.

The number of all possible integer hives with these fixed boundaries remarkably yields the corresponding LRC for the tuple.

While polynomial time algorithms exist for fixed rank Lie groups \cite{Rassart2004}, in general computing an LRC is $\# P$-complete and prompts the search for efficient estimation algorithms \cite{Narayanan2006}. Integer sampling on the hive space provides an opportunity to perform such estimations.

Numerically finding one such configuration, if it exists, is a linear programming (LP) problem on the hive polytope defined by the rhombus inequalities and the boundary data. KT have shown that finding the hive with the largest coefficients using a real valued LP can be guaranteed to give integer solutions, keeping the optimization within the HIVE space \cite{Knutson1999}. However, finding more than one solution can not be done so easily as a constrained integer programming problem. Instead, estimations of the HIVE volume prompt the need for a random walk on the HIVE space, as will be discussed in Section \ref{sec:LRC}.

\section{Hives from Hermitian Matrices}

Danilov and Koshevoy (DK), in an effort to provide a more elementary proof of the Horn problem concerning the sums of spectra of Hermitian matrix pairs, proposed a conjecture on the space of Hermitian matrices for generating real valued `hives' from an optimization process that would lend itself as an analytic continuation of the integer HIVE space \cite{Danilov2003}. This was the first explicit constructive conjecture for finding hive coefficients given two generating matrices, and they demonstrated that their conjecture holds true for select special cases: $2$ and $3$ dimensional matrices, and matrices which commute. Appleby and Whitehead (AW) later considered a modified version of this proposal for which they sketch a proof claimed to produce hives in general \cite{Appleby2014}. This would provide a novel mechanism for the exploration of hive geometries on spaces of matrix ensembles with a study of experimental observables such as the mean hive and its fluctuations, if one could find a suitable algorithm for actually implementing the optimization. We reproduce the construction to establish notation.

The optimization is concerned with finding a map from an $n\times n$ Hermitian matrix triple $(M,N,L=M+N)$ into hives with indexed entries $H_{ijk}\,,\,i=n-j-k$ in the sense of KT, where $i$ and $k$ are coordinate labels on the vertical and horizontal triangulation axes. Provided such a map exists, then the spectra $\sigma$ of the matrices $\sigma(M,N,L) = (\mu,\nu,\lambda)$ will generate a hive with the same coefficients. As in the KT case for the integer valued hives, the rhombus inequalities will again be satisfied for hive points. For completeness, we explicitly list the inequalities referenced in Eq. \ref{eq:rhombusAngles} for the right-facing, left-facing, and vertical interior rhombuses, in that order.

\begin{eqnarray}
H_{ijk} + H_{(i+1)(j-2)(k+1)} &\leq& H_{(i+1)(j-1)k} + H_{i(j-1)(k+1)} \nonumber \\
H_{ijk} + H_{(i+1)(j+1)(k-2)} &\leq& H_{(i+1)j(k-1)} + H_{i(j+1)(k-1)} \label{eq:rhombusIneq}\\
H_{ijk} + H_{(i-2)(j+1)(k+1)} &\leq& H_{(i-1)(j+1)k} + H_{(i-1)j(k+1)} \nonumber
\end{eqnarray}

Let $W$ be an $n$-dimensional complex vector space with a $p$-dimensional subspace $U$ and a $q$-dimensional subspace $V$ such that $U\subseteq V$. We denote $\mathbb{F}_{p q}(W)$ to be the set of all pairs of subspaces $(U,V)$.

We also establish a notation for orthogonal projection $\pi_U\,:\,W\rightarrow U$ such that a Hermitian operator $M$ on $W$ has an action on $U$ given by $M|_U=\pi_u \circ M \circ \pi_u$.

Provided $(M,N)$, the proposal by AW claims that the hive coefficients are generated from a maximization of traces over all possible subspaces as follows:

\begin{equation}
  H_{ijk} = \max_{(U,V)\in\mathbb{F}_{k,(j+k)}} (tr(M|_V) + tr(N|_U)) \,.
\end{equation}

We note that this departs from the conjecture by DK, wherein the original proposal optimized over pairs of \emph{mutually orthogonal} subspaces of dimensions $k$ and $j$ where $k+j \leq n$. Here, AW optimize over the set of all pairs of subspaces with the containment structure defined by $\mathbb{F}_{p q}(W)$. It is this construction that we will work with.

Although it appeared to be a constructive theorem, this optimization problem is not actually guaranteed to produce a hive in all cases. We detail below the obstruction that we have noticed, but conclude that the construction is not rendered entirely ineffective as it does work in reduced cases, such as when $M$ and $N$ are either both sorted diagonal matrices or equal to each other.

\section{Obstruction for the Subspace-Contained Optimization} \label{sec:AWCounter}

In the process of showing that the Hermitian matrix construction satisfies the rhombus inequalities, AW use the following abridged argument: by finding the subspace pairs $(U^*,V^*),(U^{**},V^{**}) $ that independently maximize each term in the left side of the inequalities in Eq. \ref{eq:rhombusIneq}, one can always form an orthonormal basis of the union $U^* \cup V^*\cup U^{**} \cup V^{**} $ such that the space can be again subdivided into subspaces that have the right dimensions for the subspace pairs (but not necessarily those that maximize the traces) $(U',V'),(U'',V'')$ of the right side of the inequalities. When the matrices $M$ and $N$ are projected into the new re-partitioned space, their traces will be at most equal to the hive coefficients of the obtuse components of the rhombus. As an explicit example presented in reference \cite{Appleby2014}, look at the right rhombus inequality:

\begin{widetext}
\begin{eqnarray}
 H_{ijk} + H_{(i+1)(j-2)(k+1)} &=& \max_{(U,V)\in F_{k,k+j}}{(tr(M|_V) + tr(N|_U))} + \max_{(U,V)\in F_{k+1,k+j-1}}{(tr(M|_V) + tr(N|_U))} \nonumber \\
 &=& [tr(M|_{V^*}) + tr(N|_{U^*})] + [tr(M|_{V^{**}}) + tr(N|_{U^{**}})] \nonumber \\
 &=& [tr(M|_{V'}) + tr(N|_{U'})] + [tr(M|_{V''}) + tr(N|_{U''})] \label{eq:AWProof} \\
 &\leq& \max_{(U,V)\in F_{k,k+j-1}}{(tr(M|_V) + tr(N|_U))} + \max_{(U,V)\in F_{k+1,k+j}}{(tr(M|_V) + tr(N|_U))} \nonumber \\
 &=&  H_{(i+1)(j-1)k} + H_{i(j-1)(k+1)} . \nonumber
\end{eqnarray}
\end{widetext}

The first and last equalities come from the hive coefficient definition; the second equality comes from the definition of the maximal subspaces; the penultimate inequality comes from the maximum principal. However, we claim that the central equality that accomplishes the re-distribution of subspaces is not always true.

The dimension counting arguments foremost do not strictly hold in the case when the subspace pairs have non-trivial overlap and the structures are generic. In this situation, it is not in general true that any reshuffling of basis elements from the first union space will guarantee that the subspaces on the right side of the inequality can be reconstructed while maintaining the total traces. We can certainly guarantee that the union of the bases for $(U^*,V^*),(U^{**},V^{**})$ is equal to that of $(U',V'),(U'',V'')$ by dimension counting. But, if in order for the sums of traces to be equivalent, it is necessary for a shared dimension to be passed from one subspace pair to another, the second subspace pair will not change rank and this may cause a hive deficiency (a failure of a rhombus inequality) as the second subspace pair will not actually have the correct dimension.

Another way to see this obstruction manifest more clearly is that it is not generically true that the rearrangement of the subspaces will result in an invariant of the sum of the traces of $M$ and $N$ upon carrying out the new projections. This is the converse of the above case, wherein we can satisfy the proper subspace dimensionality, but there is no guarantee on the total trace sums. We will construct an example for which
\begin{equation}
[tr(M|_{V^*}) + tr(N|_{U^*})] + [tr(M|_{V^{**}}) + tr(N|_{U^{**}})] > [tr(M|_{V'}) + tr(N|_{U'})] + [tr(M|_{V''}) + tr(N|_{U''})] \, .
\end{equation}
This would infringe on the crucial inequality in Eq. \ref{eq:AWProof} for the right-facing rhombus condition to hold.

Consider $(n,j,k)=(5,2,3)$, a coefficient which admits a right-facing rhombus inequality. To make our argument even more apparent, we consider positive-definite Hermitian matrices. Using a subscript notation to describe the rank of the projections with tildes denoting different subspaces in the case of dimensional equivalence, we claim that

\begin{equation}
tr(M|_5) + tr(N|_3) + tr(M|_4) + tr(N|_{4}) \stackrel{?}{>} tr(M|_{4'}) + tr(N|_{3'}) + tr(M|_{5'}) + tr(N|_{\tilde{4}'}). \label{eq:eg}
\end{equation}

Recall the variational definition of the eigenvalues of a hermitian operator acting on an $n$ dimensional vector space $W$: for spectra in weakly decreasing order $\sigma(M) = [\mu_1,\ldots,\mu_n]$ and for all $k$ such that $0 < k \leq n$,

\begin{equation}
\max_{\substack{\dim(U) = k \\  U \subseteq W}}{tr(M|_U)} = \mu_1 + \cdots + \mu_k \, .
\end{equation}

For the first trace on the left of \ref{eq:eg}, the subspace is maximal with respect to the vector space. This implies that the joint optimization between $M$ and $N$ breaks down, and we simply yield the summed spectra of $M$ with the partial spectra of $N$:

\begin{equation}
tr(M|_5) + tr(N|_3) =  \|\mu\|_1 + \nu_1 + \nu_2 + \nu_3 \\
\end{equation}

Similarly, in the second optimization on the left side of Eq. \ref{eq:eg} where the subspaces are identical, another reduction can be made directly to the eigenvalues of the summed matrices:

\begin{eqnarray}
tr(M|_4) + tr(N|_{4}) &=& tr((M+N)|_4) \\ \nonumber
 &=& \lambda_1 + \lambda_2 + \lambda_3 + \lambda_4  \,.
\end{eqnarray}

For the right side of Eq. \ref{eq:eg}, by a similar argument, we know that at most
\begin{equation}
tr(M|_{5'}) + tr(N|_{\tilde{4}'}) \leq  \|\mu\|_1 + \nu_1 + \nu_2 + \nu_3 + \nu_4 \, . \\
\end{equation}
This is due to the fact that the subspaces are not guaranteed to be those which maximize the individual traces. We substitute the maximal values as the strongest attempt for our conjectured inequality to fail.
This results in requiring that
\begin{eqnarray}
tr(M|_{4'}) + tr(N|_{3'})  &\stackrel{?}{\leq}& \lambda_1 + \lambda_2 + \lambda_3 + \lambda_4 - \nu_4 \nonumber \\
 &\stackrel{?}{\leq}& \|\mu\|_1 + \|\nu\|_1 - \lambda_5 - \nu_4   \\
 &\stackrel{?}{\leq}& (\mu_1 + \mu_2 + \mu_3 + \mu_4) + (\nu_1 + \nu_2 + \nu_3) + (\mu_5 + \nu_5 - \lambda_5) .  \nonumber \label{eq:stricterAW}
\end{eqnarray}

Equality would be necessitate that the remaining traces produce the partial eigenvalue sums of the matrices up to a difference in the smallest eigenvalue for $M$ and $N$ compared with $M+N$. It is easy to conceive of an example where this difference of the last term is $0$, so we proceed with such a case:
\begin{eqnarray}
\mu &=& (40,30,20,10,2) \nonumber \\
\nu &=& (40,30,20,10,2) \\
\lambda &=& (65,55,45,35,4) \nonumber
\end{eqnarray}

In an example like this, it would be therefore necessary for the traces to yield the partial spectra, or greater (which is not possible due to the positive definite assumption and the variational eigenvalue definition), for our claim to be false. However, it is not necessarily true that in non-maximal or non-identical subspaces that the joint optimization over the independent matrices $M$ and $N$ will yield the partial sums of their spectra over the constrained dimensions they are projected into, let alone subspaces that were only selected for their ranks as we have here. In fact, this is generically not the case. Such lower dimensional projections from the full vector space may therefore not yield equivalent traces and highly depend on the structure of the re-partitioning of the original total union space and the initial maximization.

One may ask what would happen in this example if the repartitioned subspaces \emph{were} the precise subspaces which could produce the spectral sums, giving us equality in the tightest case. We show that this does not generically hold in this example due to reasoning similarly suggested in our first argument.

In order for the exception to hold, one would need to use the same subspaces produced in the maximizations on the left side of \ref{eq:eg} for the right side (barring accidental degeneracies or special structures). That is, we would require equality of subspaces such that
\begin{eqnarray}
5 &=& 5' \nonumber \\
3 &=& 3'  \\
4 = &4'& = \tilde{4}' \nonumber \,.
\end{eqnarray}
This would be akin to the exact "redistribution" of subspaces by dimension counting referenced in \cite{Appleby2014}. However, this may not hold due to the subspace containment constraint. Necessary constraints are that
\begin{eqnarray}
3 &\subset& 5 \nonumber \\
3' &\subset& 4' \, .
\end{eqnarray}
Yet, due to our repartition, we would require an additional constraint that is \emph{not} enforced:
\begin{equation}
  3 \subset 4 \, .
\end{equation}
Thus, in order to fix equality of traces, we must also have a guarantee that the subspace which maximizes $tr(M+N)$ also maximizes $tr(N)$, and this is not generically true. For matrices which do not have such an additional structure, our inequality in Eq. \ref{eq:eg} is true, which violates the original equality used in Eq. \ref{eq:AWProof} by AW to support their theorem.

This condition leads us to new matrix domains where the optimization \emph{is} guaranteed to hold: matrices with the same eigenstructure ranking. Sorted diagonal matrices, or matrices for which $M=N$, are two examples of classes of matrices wherein this constraint will be automatically satisfied. Although these are more trivial classes than one would hope to study in general, it will still give some avenue to utilizing the optimization problem presented above in order to study hive structures.

\section{Optimization Space Transformation for a Practical AW Implementation}

Sampling over all such projections into the constrained subspaces as per the AW construction would be akin to an optimization problem on a flag manifold, or a collection of ordered sets of vector subspaces. This is a difficult process to accomplish algorithmically, and it could be considered a generalization of independent subspace analysis \cite{Nishimori2006}. Instead, we introduce a map of the AW optimization problem into a numerically equivalent problem on a product space of Grassmannians, on which one can perform a Riemannian optimization scheme. We produce the analytic Euclidean gradients and Hessians necessary to implement a traditional gradient-descent or trust-region algorithm, which we can then transform to work on the product space and return the individual hive coefficients.

We begin by trivially turning our maximization into a minimization for the following algorithm:
\begin{eqnarray}
  H_{ijk} &=& \min_{(U,V)\in\mathbb{F}_{k,(i+k)}} -(tr(M|_V) + tr(N|_U)) \nonumber \\
  &=& \min_{(U,V)\in\mathbb{F}_{k,(i+k)}} -(tr(\pi_V M \pi_V) + tr(\pi_U N \pi_U))
\end{eqnarray}

We introduce an explicit form for the projection with $\vec{\alpha}_i$ defined implicitly as follows:
\begin{eqnarray}
  \pi_a &=& A(A^TA)^{-1}A^T \\
  &\text{s.t.}& \nonumber \\
  A &=& \text{col}(\vec{\alpha}_i)\,|\,i\in\{0,\ldots,\text{rank}(a)\} \nonumber \\
  a &=& \text{span}(\vec{\alpha}_i) \nonumber \, .
\end{eqnarray}

Expanding out, we can now optimize over new variables:

\begin{equation}
  H_{ijk}= \min_{\substack{A=\text{col}(\vec{\alpha}_i) \,|\,V=\text{span}(\vec{\alpha}_i) \\ B=\text{col}(\vec{\beta}_i) \,|\,U=\text{span}(\vec{\beta}_i) \subseteq V}}
  -(tr(A(A^TA)^{-1}A^T M A(A^TA)^{-1}A^T) + tr(B(B^TB)^{-1}B^T N B(B^TB)^{-1}B^T)) \, .
\end{equation}

Now we decompose $A$ into the subspace containing $U$ and any orthogonal space, $A=\rv{B}{\tilde{A}}$:

\begin{eqnarray}
  H_{ijk} = \min_{\substack{\tilde{A}=\text{col}(\vec{\alpha}_i) \,|\,\text{span}(\vec{\alpha}_i) \perp U \\ B=\text{col}(\vec{\beta}_i) \,|\,U=\text{span}(\vec{\beta}_i) \subseteq V}}
  &-&(tr(\rv{B}{\tilde{A}}(\rv{B}{\tilde{A}}^T\rv{B}{\tilde{A}})^{-1}\rv{B}{\tilde{A}}^T M \rv{B}{\tilde{A}}(\rv{B}{\tilde{A}}^T\rv{B}{\tilde{A}})^{-1}\rv{B}{\tilde{A}}^T) \nonumber \\
  &+& tr(B(B^TB)^{-1}B^T N B(B^TB)^{-1}B^T)) \nonumber \\
  &\equiv& \min_{\substack{\tilde{A}=\text{col}(\vec{\alpha}_i) \,|\,\text{span}(\vec{\alpha}_i) \perp U \\ B=\text{col}(\vec{\beta}_i) \,|\,U=\text{span}(\vec{\beta}_i) \subseteq V}} -(f(B,\tilde{A},M)+g(B,N))\, ,
\end{eqnarray}
for functions $f$ and $g$ defined as the respective portions of the cost function.

It is at this step that we notice that if one were to drop the orthogonality constraint, the optimization over $\tilde{A}$ and $B$ would live on a product space of Grassmannians,

\begin{equation}
  (\tilde{A},B) \in \mathfrak{Gr}_i^n \times \mathfrak{Gr}_k^n \, .
\end{equation}
Formally, we can impose the orthogonality by introducing a cost penalty $C$ of the following form:

\begin{equation}
  C = \lim_{\xi\rightarrow \infty} e^{\xi ((i+k)-\text{rank}(A))} \,.
\end{equation}

This gives an infinite penalty whenever $A$ is not full rank, or in other terms, whenever $\tilde{A}$ has a span that becomes subspace degenerate with the column span of $B$. This cost is not possible to implement or account for numerically. However, if one were randomly sampling matrices from the product space with double precision numerics, this degeneracy would occur with vanishing probability and $(A^TA)^{-1}$ will never be singular. Yet, in an optimization scheme, one has to worry about whether the optimal flow will drive toward such a degeneracy. A simple heuristic argument can convince the reader that this should not happen for symmetric positive-definite matrices, however. The projection will yield a positive semi-definite matrix with null columns on orthogonal projected subspaces. Larger rank subspaces (i.e., whenever $\tilde{A}$ is full) contribute more to the maximizing trace we are optimizing. Therefore, the minimum (negative trace) should be driven \emph{away} from these degenerate configurations and we should not expect numerical instability arising from nearly singular matrices. In practice, this is confirmed.

In this way, we have transformed the problem from constrained optimization on a flag manifold into an unconstrained optimization over a product of Grassmannians.

\subsection{Testing the Grassmannian Map}

  Once we have computed the Euclidean gradients and Hessians, the \emph{MANOPT} matlab package \cite{Boumal2014} can internally perform the required transformation to relate our operators with the operators on the product manifold of Grassmanians.

  We have checked the accuracy of our expressions against numerically calculated gradients and Hessians. Our expression is confirmed to stay within the correct tangent space and have the proper magnitude.

\begin{figure}[H]
  \centerline{\includegraphics[scale = .5]{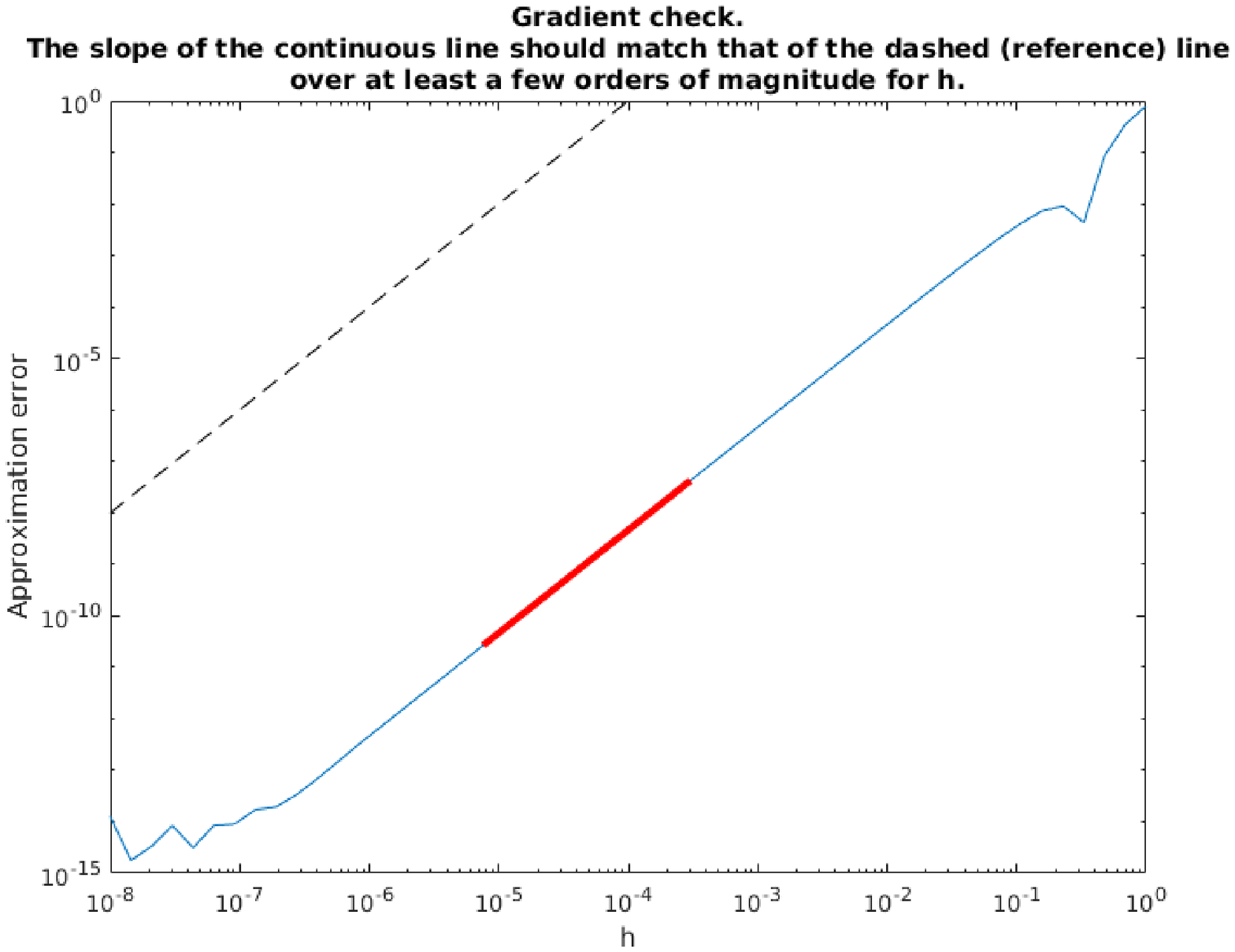}\includegraphics[scale = .5]{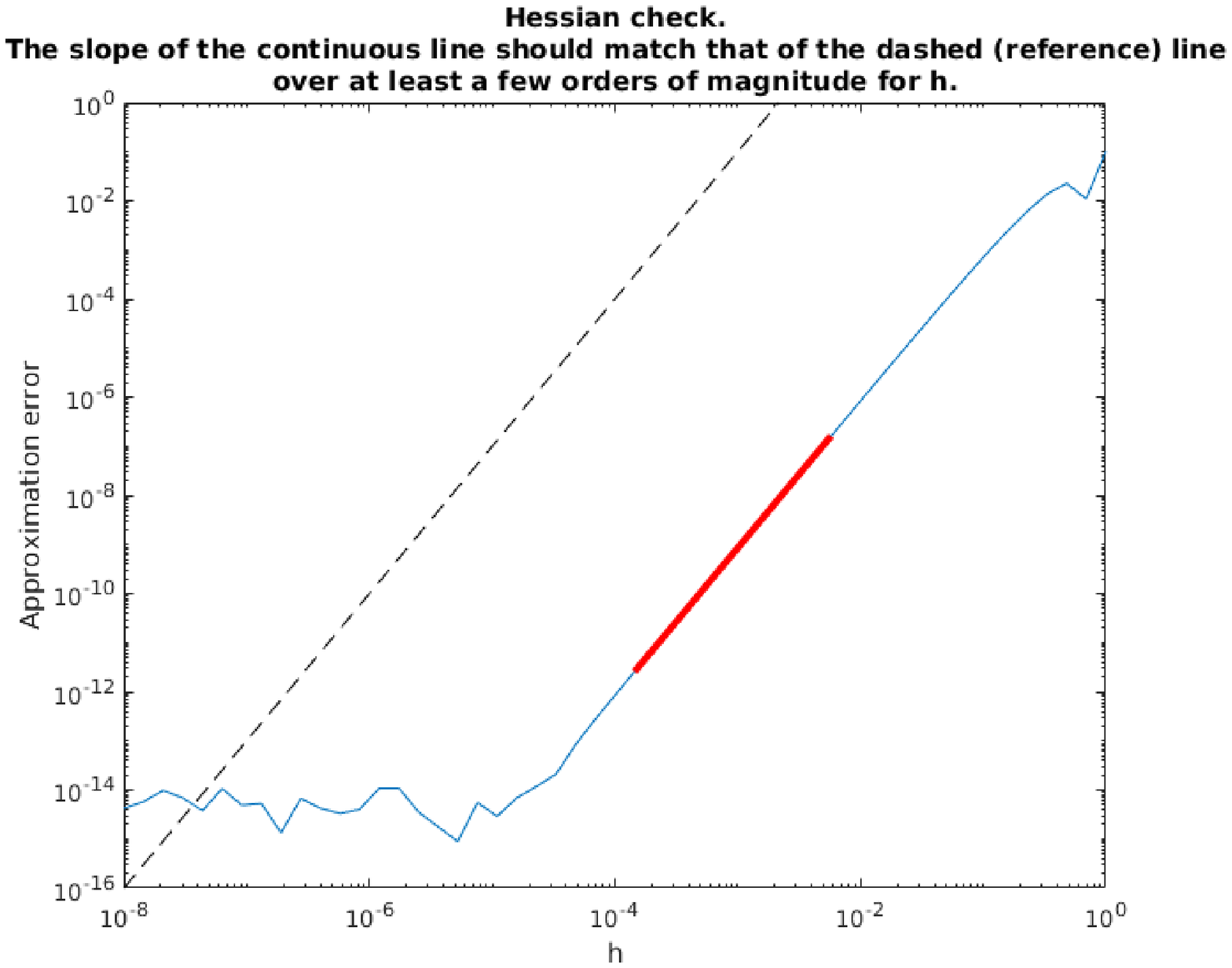}}
  \caption{Gradient and Hessian accuracy test, indicating strong agreement between slow and potentially unstable numerical approximations and our derived results.}\label{fig:Numerics}
\end{figure}

  The graphical and numerical output from \emph{MANOPT} for a sample optimization is reproduced below and illustrated in Fig \ref{fig:Numerics}.
  \begin{quote}[Gradient]
    The slope should be 2. It appears to be: 1.99989. If it is far from 2, then directional derivatives might be erroneous.
    The residual should be 0, or very close. Residual: 4.51188e-16. If it is far from 0, then the gradient is not in the tangent space.
\end{quote}
\begin{quote}[Hessian]
The slope should be 3. It appears to be: 2.99794. If it is far from 3, then directional derivatives or the Hessian might be erroneous.
The residual should be zero, or very close. Residual: 9.48494e-16. If it is far from 0, then the Hessian is not in the tangent plane.
<d1, H[d2]> - <H[d1], d2> should be zero, or very close. 	Value: 0.541143 - 0.541143 = -9.76996e-15. If it is far from 0, then the Hessian is not symmetric.
\end{quote}

\begin{figure}[H]
  \centerline{\includegraphics[scale = .5]{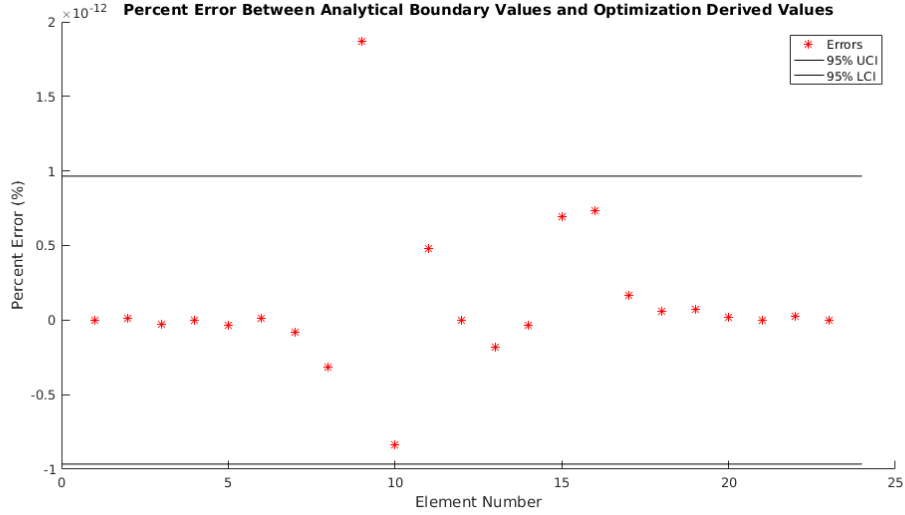}}
  \caption{A hive boundary accuracy test for an 8-D vector space, illustrating that our optimization can achieve the correct hive boundary values on the order of machine precision percent error. Upper and lower 95\% confidence intervals are illustrated.}\label{fig:optbndaccuracy}
\end{figure}

A further check that our analytical work is producing the correct hive coefficients can be performed by looking at the boundary elements of the hive. These elements are known by construction using simple algebra after the spectrum of the matrices are determined; however, we can ask our algorithm to compute these elements by way of the optimization scheme instead. Finding the boundary elements of the hive involves looking at the special cases when the contained subspace $U$ is empty, when it is the same as $V$, and when there is no dimensional reduction projection for the larger of the subspaces and $V$ spans $W$. Agreement with our optimization algorithm and the boundary data would therefore imply that our prescription is not only valid on the special cases when the Grassmann manifold dimensions $i$ and $k$ are $0$, but in a general case whenever the product manifold is full rank $i+k=n$ and is thus a very good indicator. We find excellent agreement between the two, as illustrated in Fig. \ref{fig:optbndaccuracy}.

\section{Statistics on Generating Hives}\label{sec:hivePics}

\begin{figure}[H]
  \centerline{\includegraphics[scale = .6]{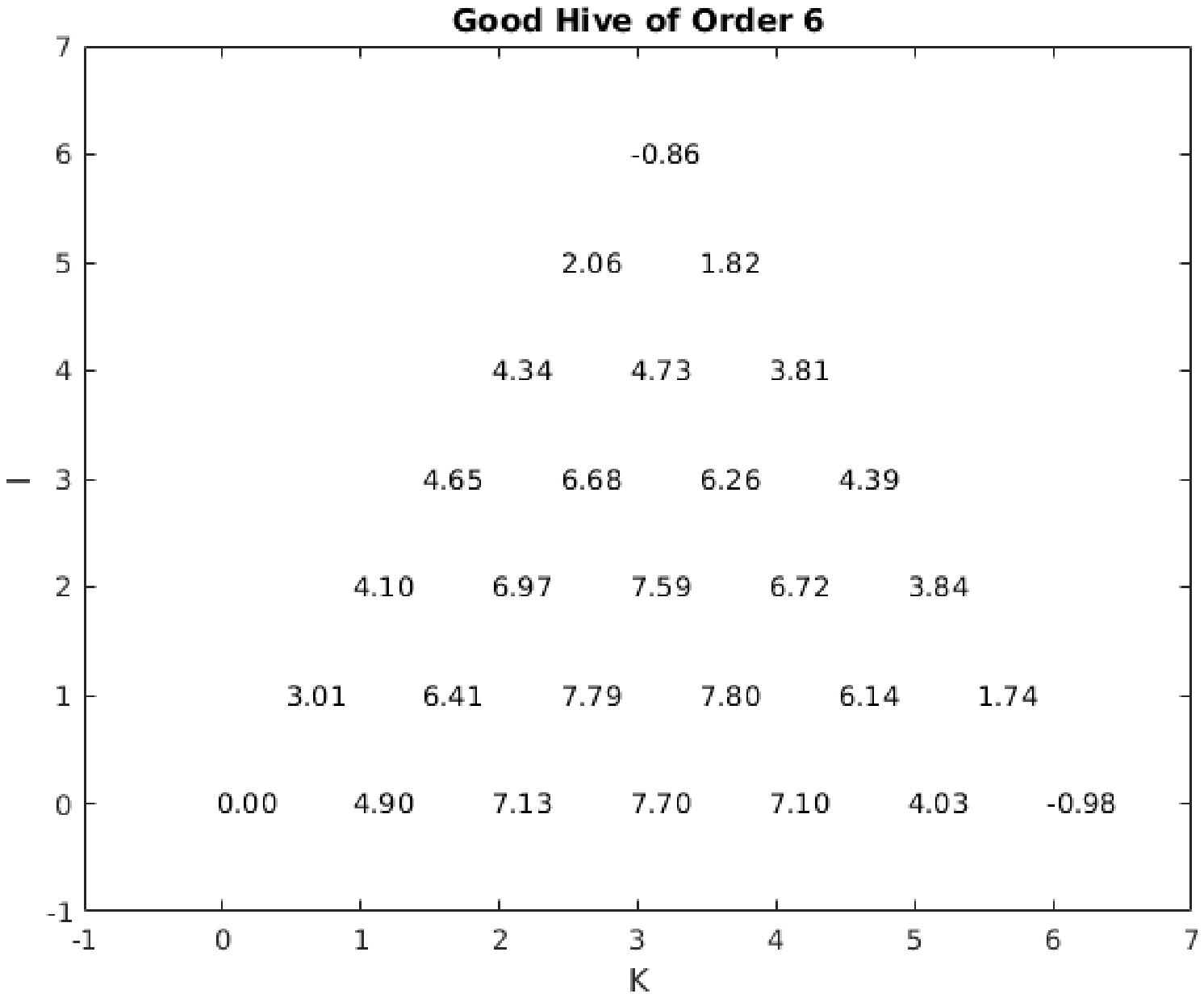} \includegraphics[scale = .6]{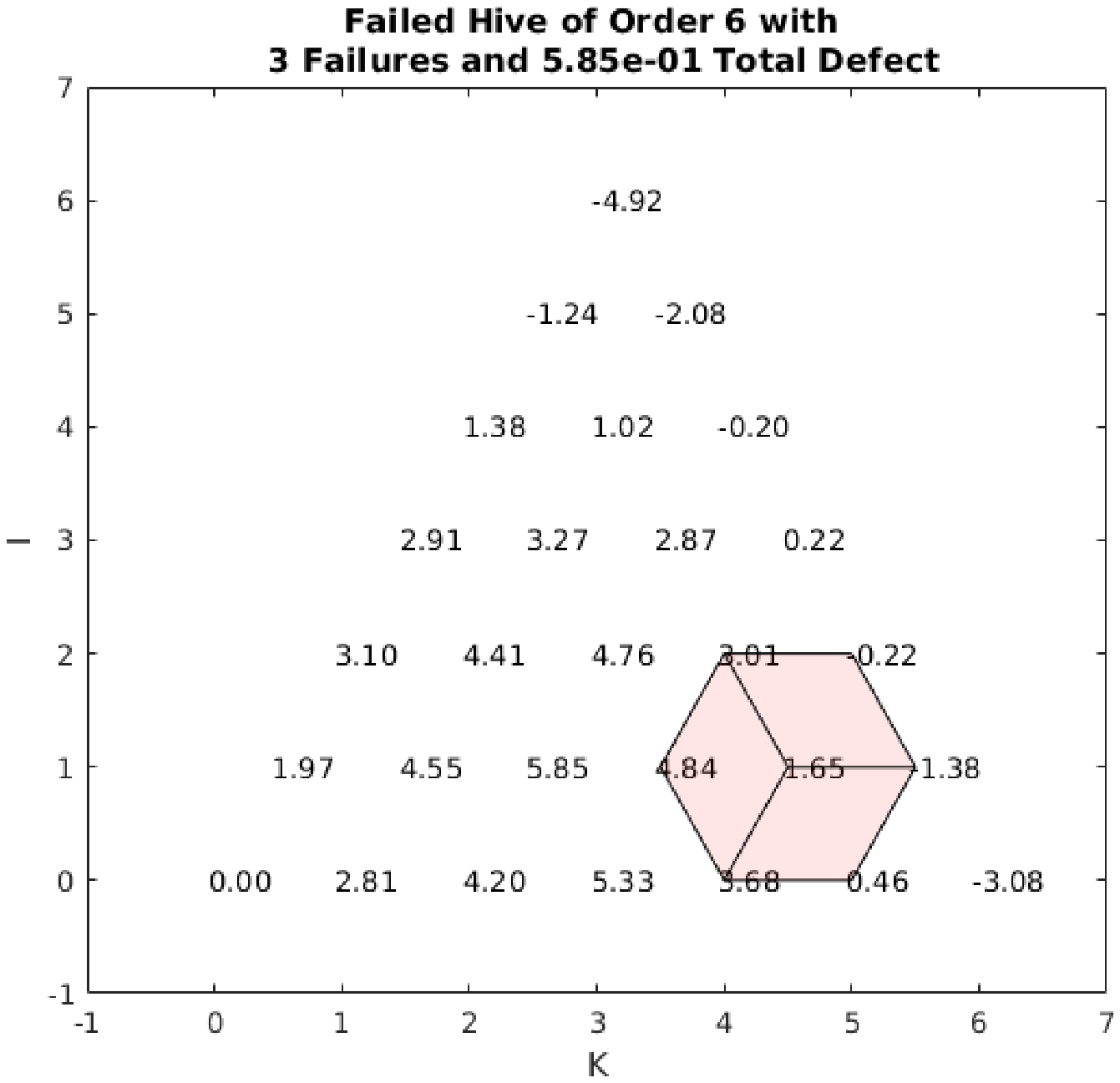}}
  \caption{On the left we show a proper hive generated by our optimization map using matrices from $GOE(6)$; on the right is an improper hive from the same ensemble with the rhombus failures highlighted. This failure was unable to be rectified by re-optimizations, suggesting a genuine numeric counterexample to the AW proposal.}\label{fig:GOE6}
\end{figure}

We have confirmed numerically that using random Hermitian matrices does not always guarantee a hive, in accordance with our theoretical analysis. Examples of a success and failure from the Gaussian Orthogonal Ensemble $GOE(6)$ are shown in Fig \ref{fig:GOE6}.

\begin{figure}[H]
  \centerline{\includegraphics[scale = .6]{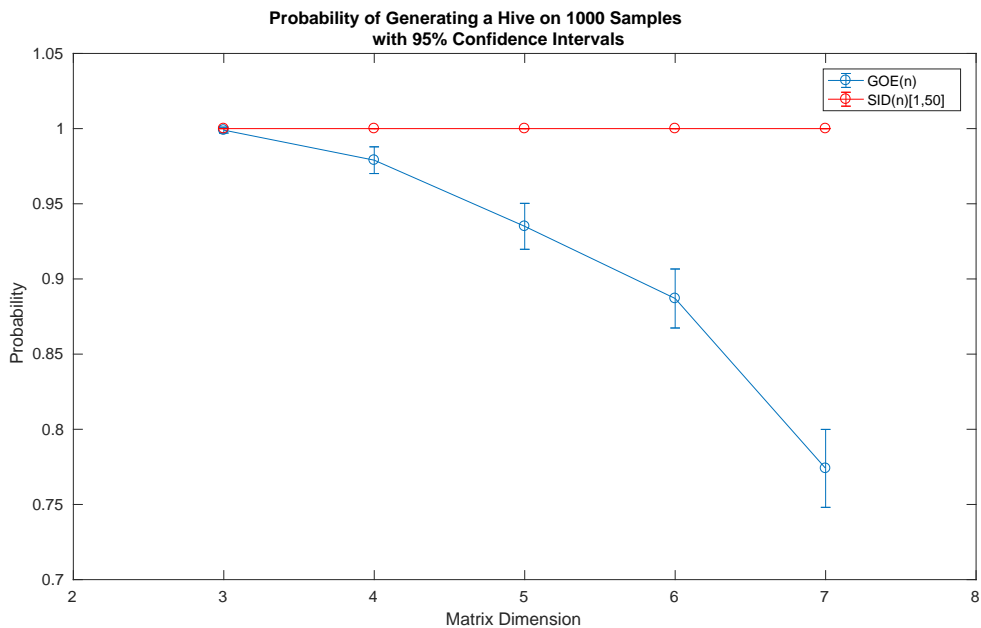}}
  \caption{Probability and 95\% confidence intervals for finding a hive on independent $GOE(n)$ matrices, and on sorted integer diagonal matrices with entries within $[1,50]$ ($SID(n)[1,50]$)}\label{fig:Probs}
\end{figure}

To get a sense for the scaling of the ability to find a hive as a function of matrix dimension in the case when the construction is not guaranteed, we sampled over such a set. We allowed for up to $5$ re-optimizations per boundary in an attempt to eliminate accidental instabilities or local minima which were not hives. Additionally, each hive coefficient allowed for 5 re-optimizations in order to handle any individual numerical instabilities. If no hive could be found after such re-optimization attempts, the matrix pair was counted as not being able to produce a hive with respect to this optimization algorithm. The results are shown in Fig \ref{fig:Probs}, showing a stark decline in probability for $GOE(n)$ with increasing matrix dimension. However, if we let $M=N$ or consider sorted diagonal matrices (also illustrated in Fig \ref{fig:Probs}) the algorithm appears to work with almost certain probability.

As mentioned in Sec.\ref{sec:AWCounter}, one way to guarantee that the resulting traces are invariant after redistributing the subspace pair is to let $M=N$. The equality of the matrices is a structure that prevents such a failure as the traces are now maximized over the same subspaces, just with a dimensional restriction. We will always be able to trade dimensions to construct the new subspaces while maintaining the traces as the spaces overlap.

If we use this setup numerically, for convergent optimizations, we yield hives with extremely strong probability. Numerical instabilities in the high dimensional optimizations may induce hive failures, but this is due to non-convergent optimizations discussed in Sec.\ref{sec:numericalInstab} and not failures of the maximization condition. Of course, the optimization is in theory far easier now and the form \emph{could} be reduced significantly when we exploit $M=N$. However, in our implementation, the optimization has not been adjusted to exploit the symmetry and it is still optimizing over the full space assuming $M$ and $N (=M)$ are independent. This is further strong evidence for the obstruction illuminated in Sec.\ref{sec:AWCounter} and that the equality of the traces in the new subspaces is necessary for a hive to form.

This setup is restrictive, however, as it forces the main boundaries to have identical eigenvectors. One way to satisfy the new trace requirements while allowing for additional freedom is to simply use sorted diagonal matrices. Let the desired eigenvalue distributions be specified in a weakly ordered form, and let the matrices $M$ and $N$ be the diagonal matrices built from the eigenvalue distributions.
\begin{eqnarray}
 M &=& diag([\mu_1 \cdots \mu_n]) \nonumber \\
 N &=& diag([\nu_1 \cdots \nu_n])
\end{eqnarray}

We note that for the integer valued case, the third boundary will also naturally have integer values due to this structure. It is otherwise not generically true that the sum of two symmetric matrices with integer spectra yields a third matrix with integer spectra.
The ordering does not have to be descending--it is only important that the matrices have coherent ordering. This is enough structure to ensure that the trace equality will be met, as subspace containment during the first optimization (similar to the identical matrix case) will be guaranteed due to the ordering. In discussing an application to honeycombs and crystals, AW map a tuple of eigenvectors into sorted diagonal matrices for their example. Their example would still hold in light of our findings, but the prescription would not generate a hive if they selected any other matrix representation that happened to have the same spectra but not be diagonal (or otherwise have the correct ordered relationship between the subspaces and eigenstructures) \cite{Appleby2014}. This setup was again confirmed numerically as shown in Fig. \ref{fig:Probs}.

\begin{figure}[H]
  \centerline{\includegraphics[scale = .5]{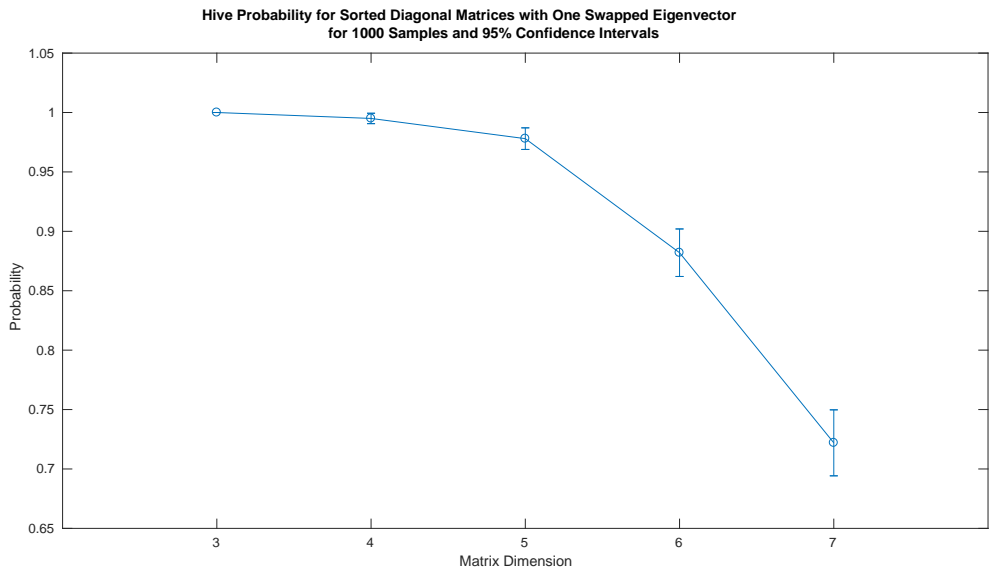}}
  \caption{Probability and 95\% confidence intervals for finding a hive on sorted integer diagonal matrices with entries within $[1,50]$, where one of the matrices had its highest and lowest eigenvalues flipped.} \label{fig:unsorted}
\end{figure}

As an additional test, if we reorder the highest and lowest eigenvalues in just one of the matrices as follows:
\begin{equation}
M \rightarrow diag([\mu_n \,\, \mu_2 \cdots \mu_{n-1} \,\, \mu_1]) ,\nonumber
\end{equation}
the optimization will fail to guarantee to produce hives as the subspaces that maximize the first hive coefficients no longer have the right containment structure for a proof of the form suggested above. Fig. \ref{fig:unsorted} illustrates the decline in probability.

In the original conjecture by DK, commuting matrices were proven to be a domain over which their proposal holds true \cite{Danilov2003}. Our numerical results with unordered diagonal matrices suggest that the AW construction does not hold for commuting matrices, making the subspace-contained projection condition less applicable than the DK proposal. We can make no statements about the DK conjecture based on our results, as the optimization spaces are mutually exclusive.

We do mention, however, that DK proved their conjecture for matrices of dimension 3. Our numerical results appear to confirm that the AW construction also holds in general for $3 \times 3$ matrices, where there is only a single interior hive point.

Thus, over this constrained space of matrices, we can still probe the geometric and probabilistic structure of some matrix ensembles and produce hives using the AW prescription.

\section{Optimization Instabilities}\label{sec:numericalInstab}
\begin{figure}[H]
  \centerline{\includegraphics[scale = .6]{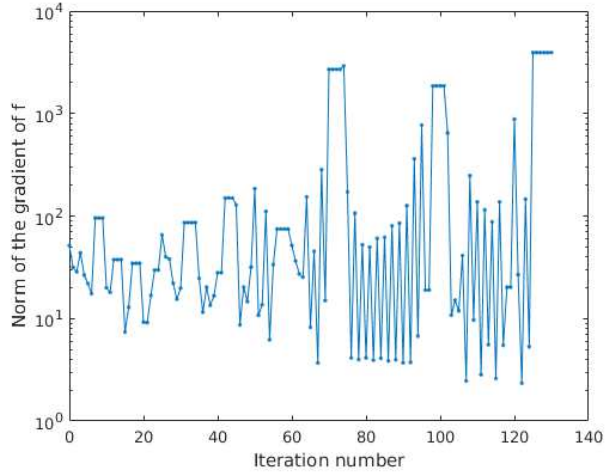}}
  \caption{A failure of convergence in the gradient norm from the trust-region algorithm iteration steps}\label{fig:gradDiv}
\end{figure}

\begin{figure}[H]
  \centerline{\includegraphics[scale = .5]{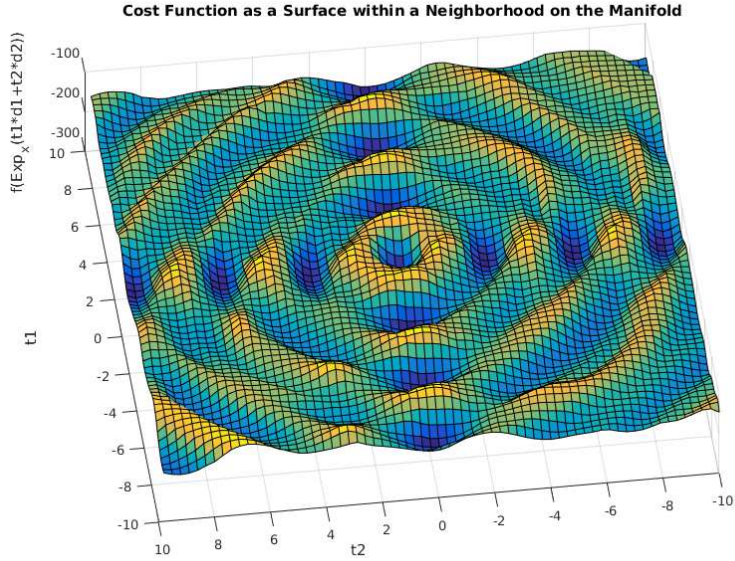}}
  \caption{An example cost function 2-D subspace in a local neighborhood about the optimum}\label{fig:costFctnSurface}
\end{figure}

We note that in terms of numerical instabilities, we occasionally see an oscillation in the gradient norm followed by a fixed-point that is accompanied by a divergence of the joint-matrix inverse. This is illustrated in Fig \ref{fig:gradDiv}. This would be indicative of the type of divergence described above when considering the missing orthogonality constraint on the product manifolds, but given the same matrix data for positive-definite matrices, re-optimizations can resolve this divergence and produce a valid hive. This indicates that the divergence is, if fact, not located at the optimum, but is instead due to unstable attractors or simply overshoot in the algorithm. As can be seen in Fig \ref{fig:costFctnSurface}, the symmetries of the Grassmannians as well as the properties of the traces in the cost function induce symmetries in the cost function space with a variety of nearby minima. Oscillations between such minima are possible, as well as issues handling the multi-dimensional saddles that clearly arise. We claim that these are the driving factors behind any instabilities in the optimization, and not any null-space issues arising from the formulation map itself (at least for the matrix ensembles with positive spectra). Finding a description of the optimization problem that does not contain redundant degrees of freedom may ameliorate some of these issues.

\section{An Open Study of Hive Properties on Accessible Matrix Ensembles}
Lastly, we examine a few matrix ensembles that give hives with extremely high probability for their average hive surface, average Gaussian curvature, and average mean curvatures. These are a few natural observables to compute when studying the random surfaces generated by the hive ensembles. We set a standardized placement of the hive points on a $2D$ plane, and triangulate the interiors of the numerical hives using a Delaunay triangulation with the hive values as heights in an orthogonal direction. We then randomly sample from the matrix ensembles until we gather a fixed number of good hive samples. We render the algebraic average of the hive surfaces, as well as the mean curvatures. The curvatures on the edges are set to $0$ by default, creating edge effects that should be noted.

\subsection{$GOE(n)$}
Pairs of identical matrices from the GOE create a highly symmetric parabolic surface. As we expect the distribution of eigenvalues to follow the Wigner semicircle law in the large matrix dimension limit, the cumulative sums on the angled boundaries also exhibit this behavior. With the bottom left corner normalized to 0, and with an even average of signed eigenvalues from the Wigner law, the tip of the triangulation also averages to 0 and the projected boundary curves resemble that of a parabola provided the right triangulation spacing \cite{Wigner1955}.
\begin{figure}[H]
  \centerline{\includegraphics[scale = .6]{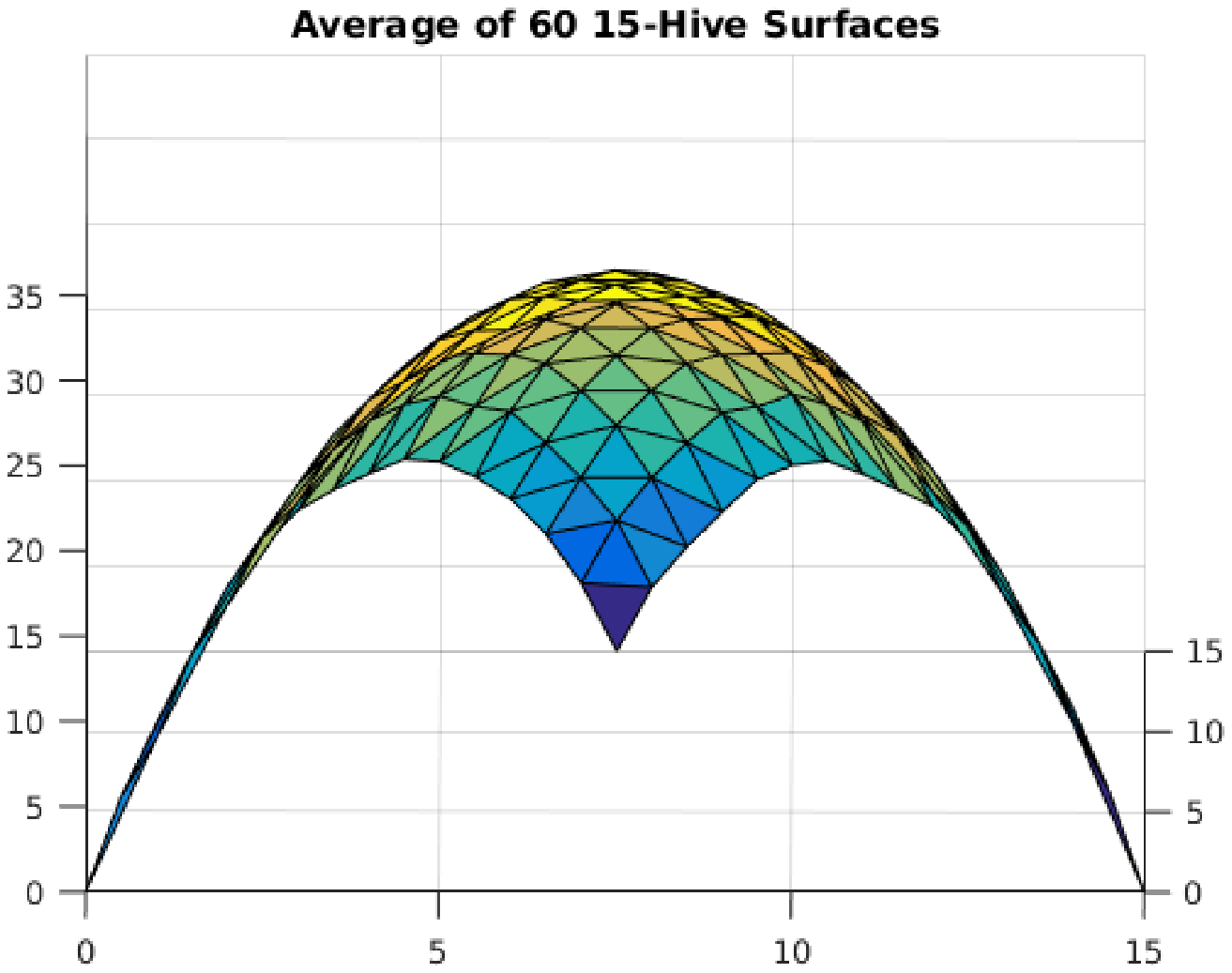} }
  \caption{Average hive surface on 60 identical pairs of random GOE(15) matrices, illustrating the `bowl-like' nature of the hives}
  \centerline{\includegraphics[scale = .6]{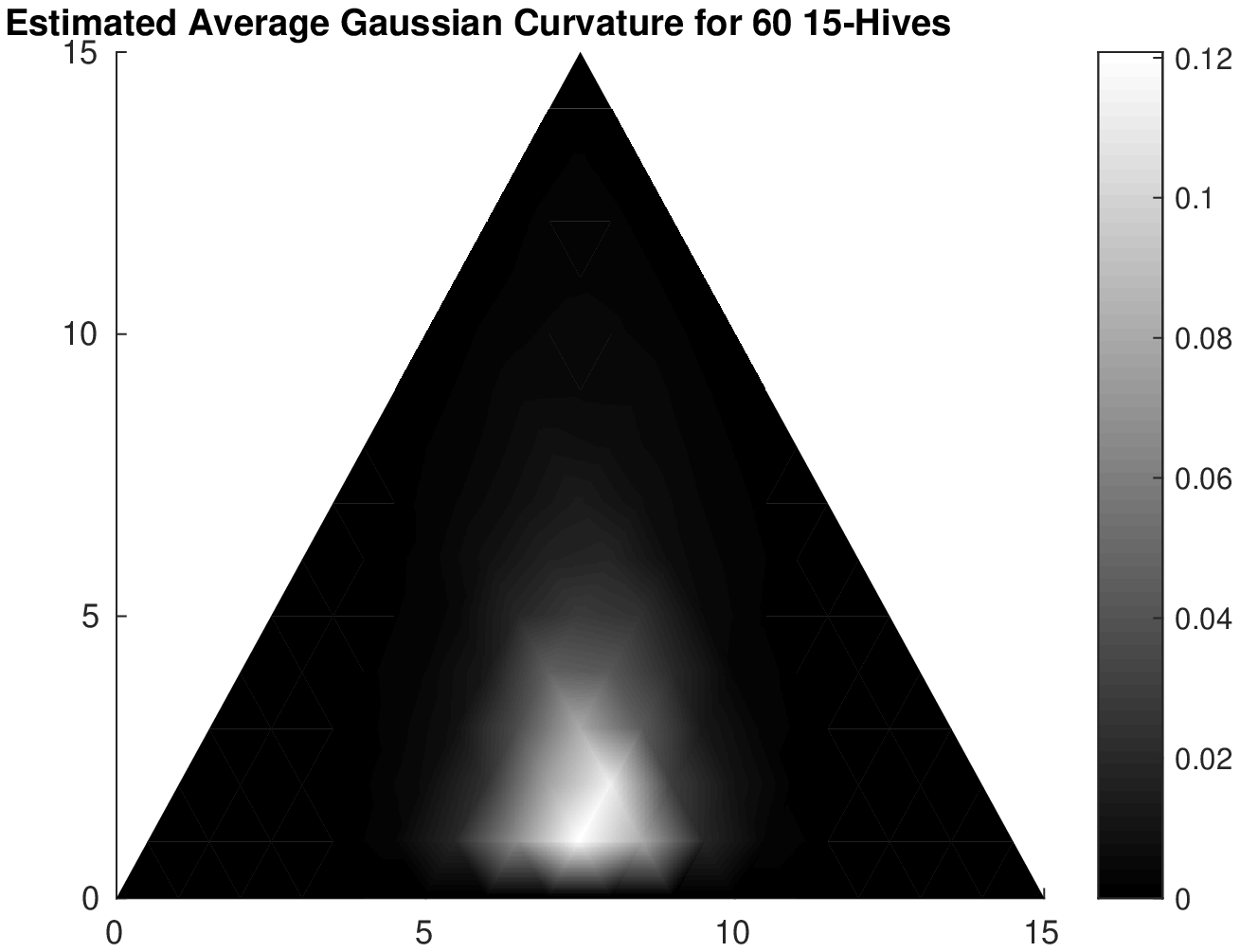} \includegraphics[scale = .6]{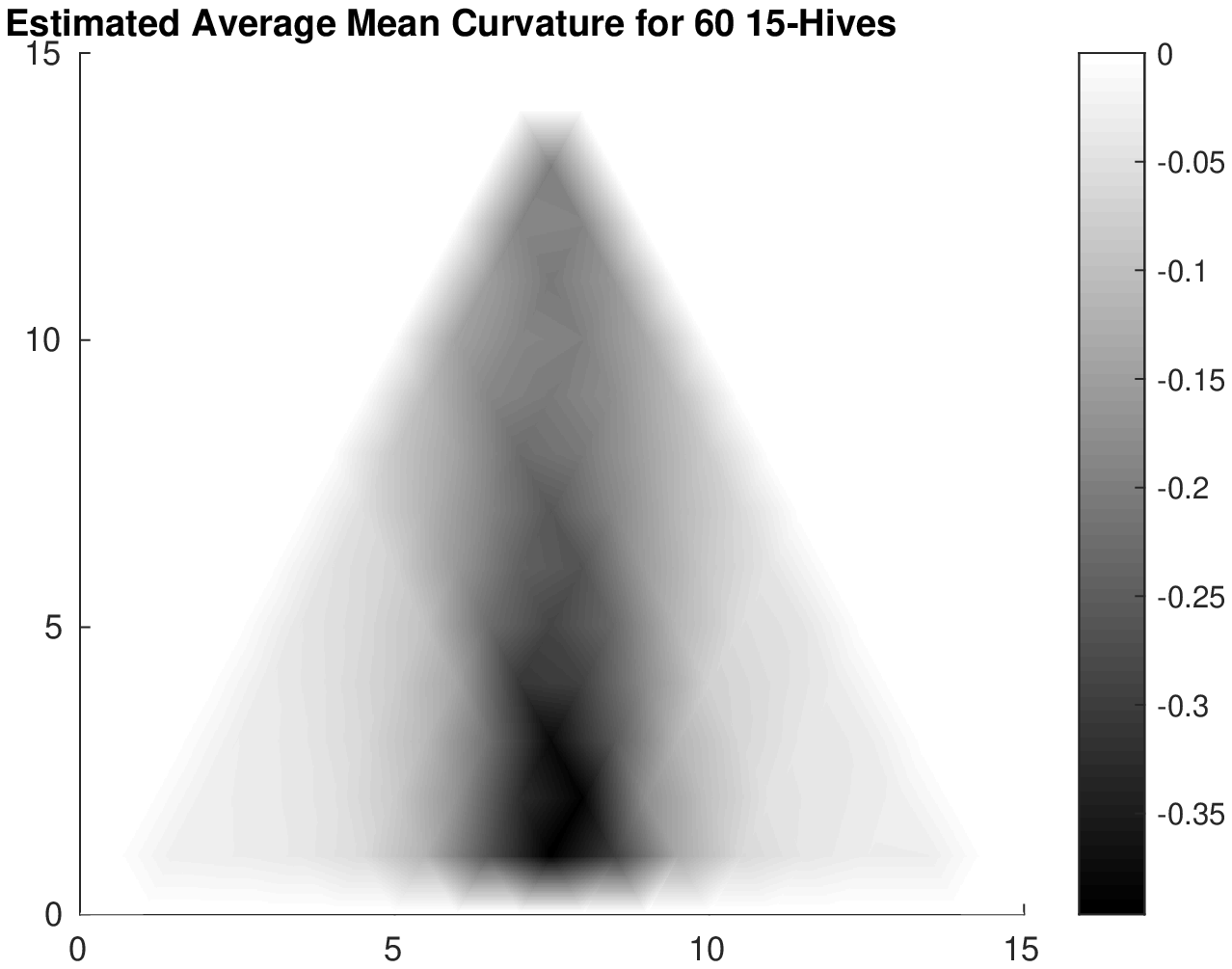}}
  \caption{Average gaussian and mean curvatures on 60 identical pairs of random GOE(15) matrices}
\end{figure}

\subsection{Integer Diagonal}

Without the oppositely signed eigenvalues, the diagonal matrices produce a surface of nearly constant $0$ Gaussian curvature and a tapering mean curvature which becomes nearly flat at its highest point.
\begin{figure}[H]
  \centerline{\includegraphics[scale = .6]{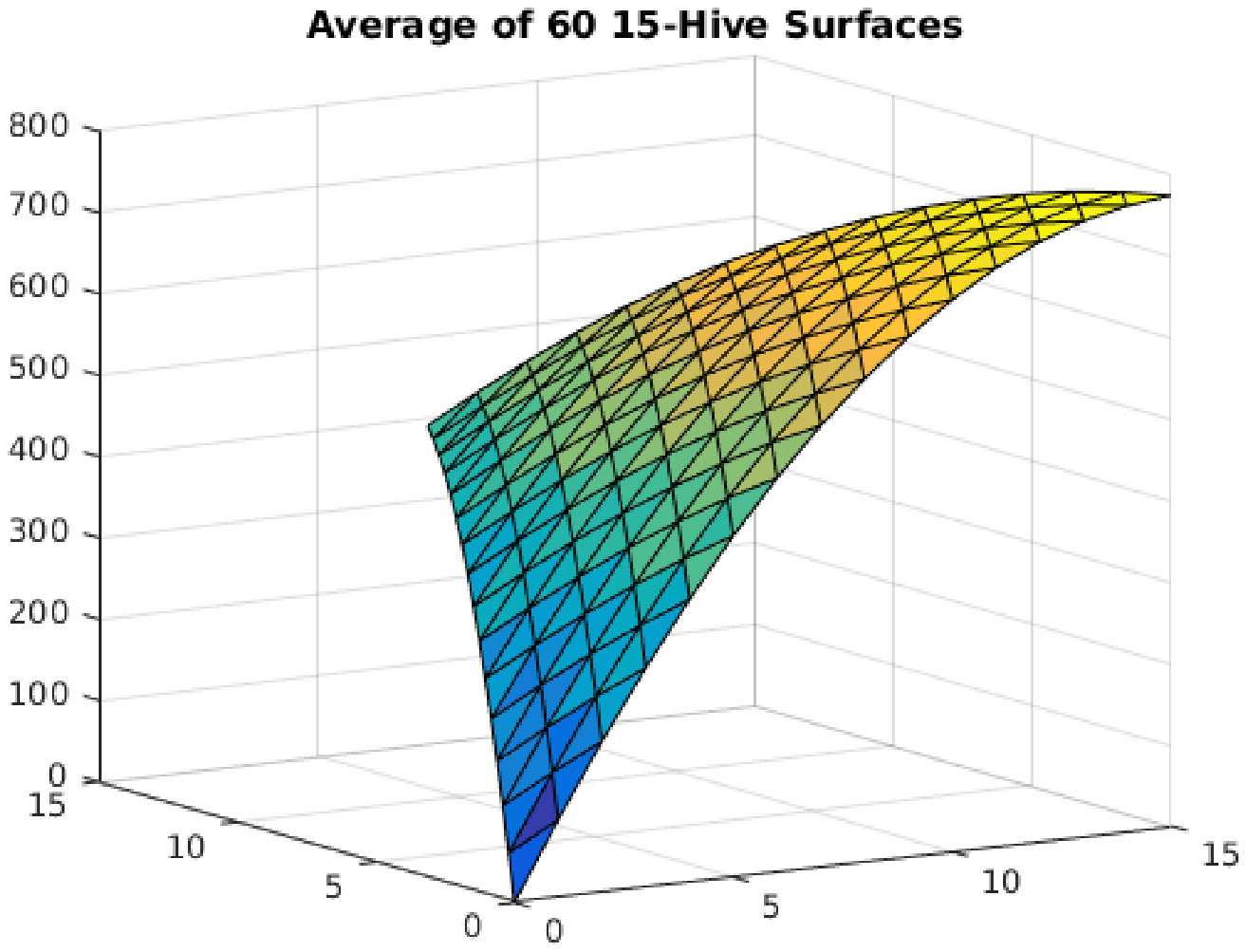}}
  \caption{Average hive surface on 60 pairs of weakly decreasing positive integer $15\times15$ diagonal matrices, uniformly distributed on $[1,50]$}
  \centerline{\includegraphics[scale = .6]{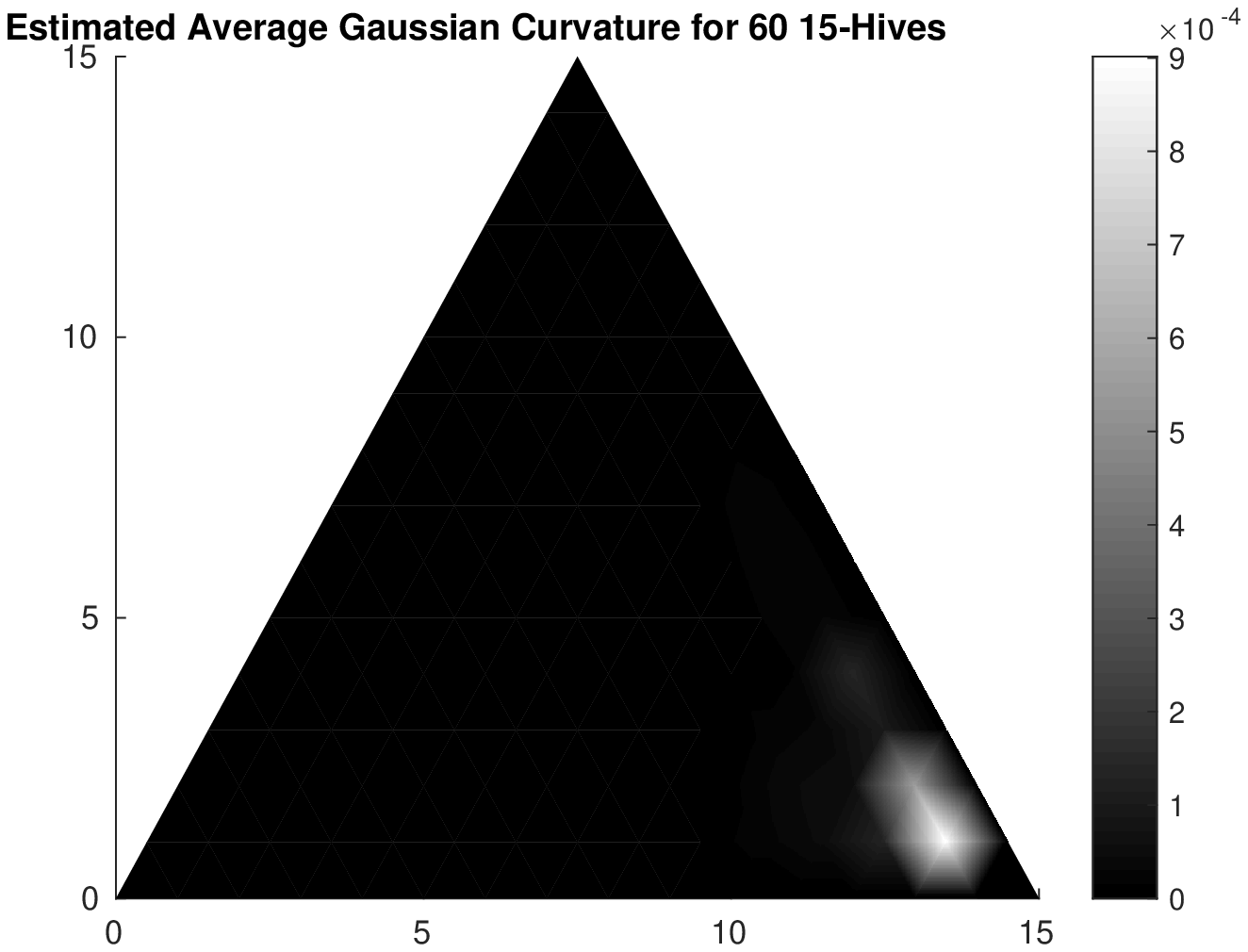} \includegraphics[scale = .6]{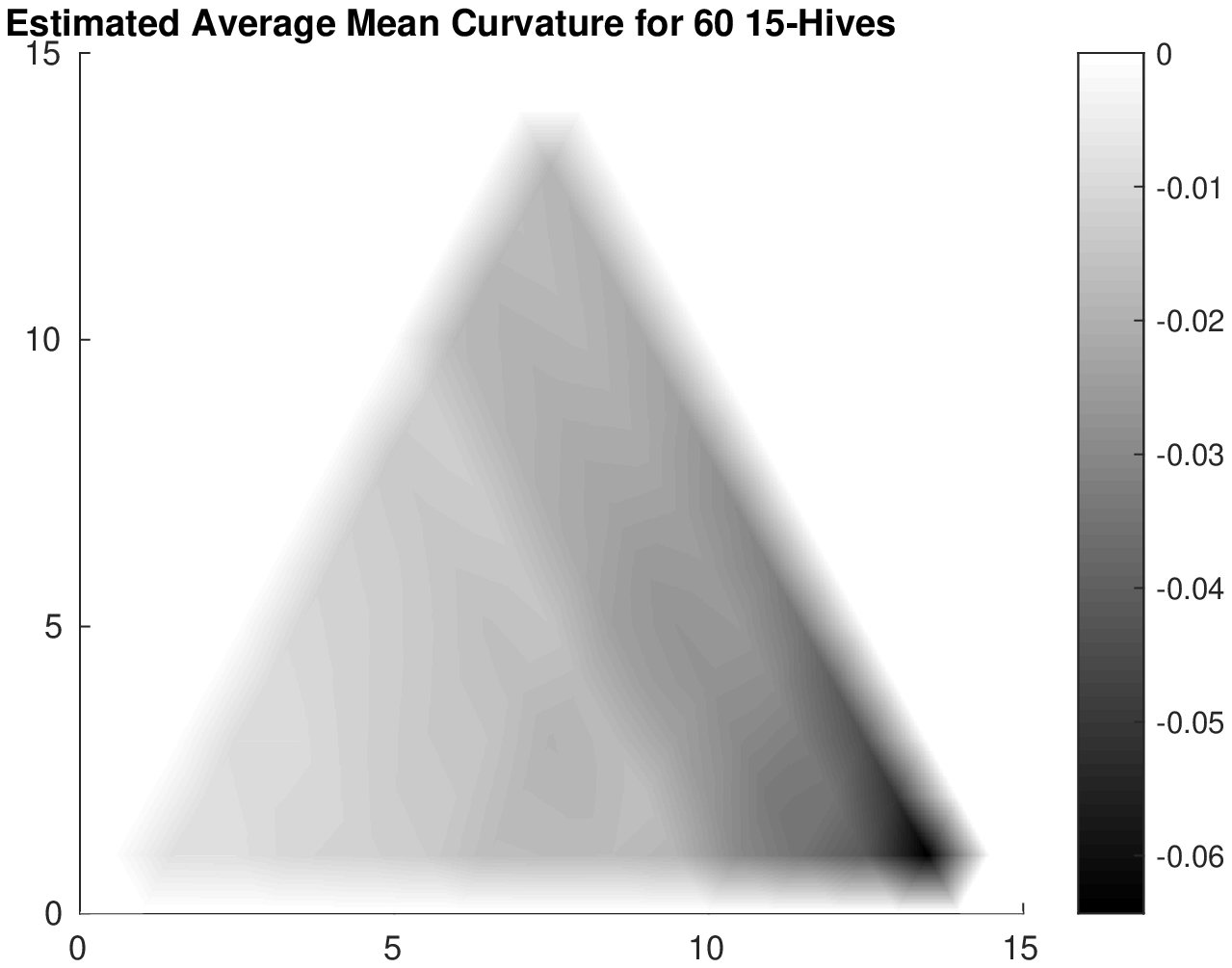}}
  \caption{Average gaussian and mean curvatures on 60 pairs of weakly decreasing positive integer $15\times15$ diagonal matrices, uniformly distributed on $[1,50]$}
\end{figure}

\subsection{$SPD(n)$}

The normally distributed symmetric positive definite (SPD) ensembles create a significantly more warped surface than the diagonal matrices, but in contrast, the diagonally dominant SPD matrices show a nearly flat surface over the whole triangulation.

\begin{figure}[H]
  \centerline{\includegraphics[scale = .6]{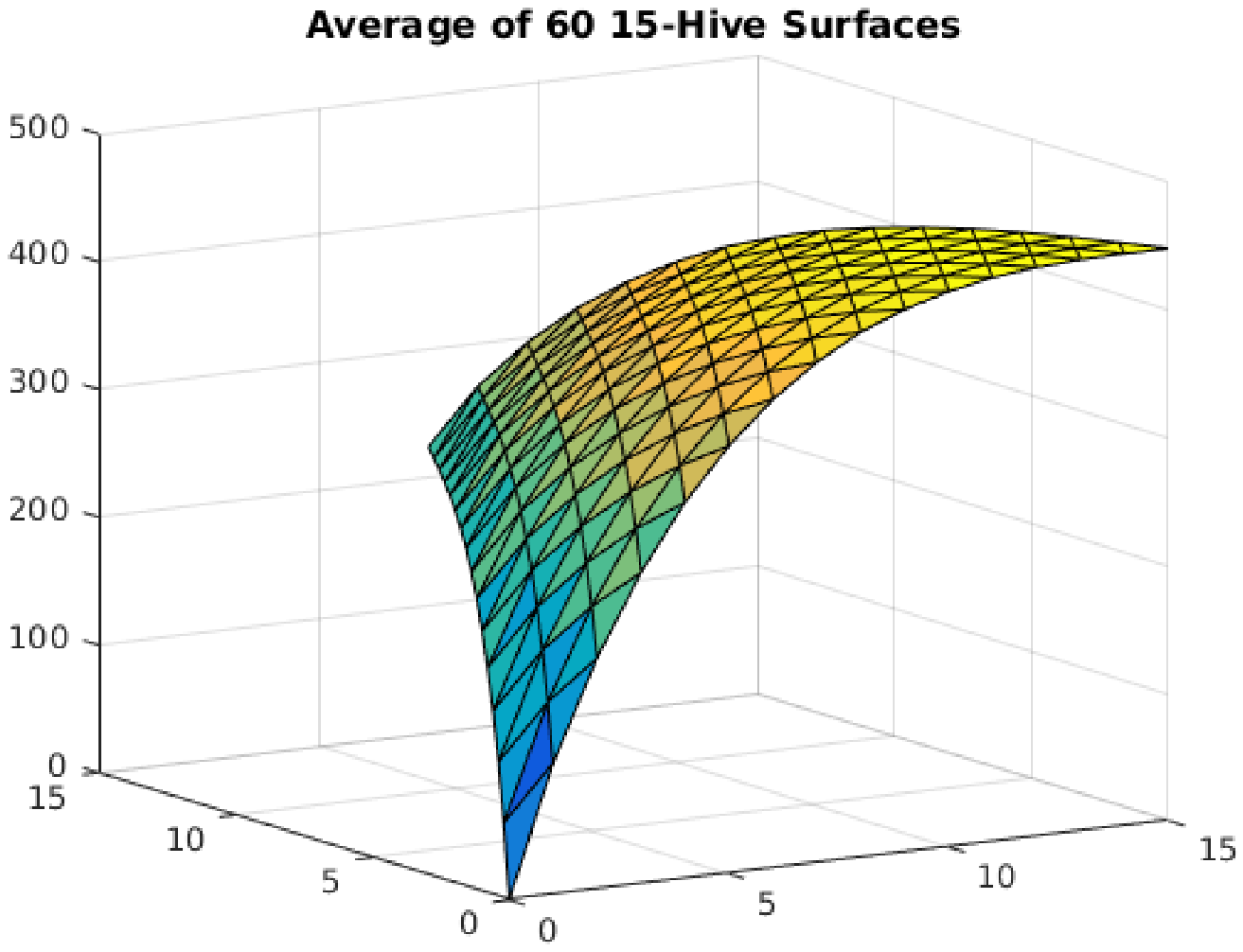} \includegraphics[scale = .6]{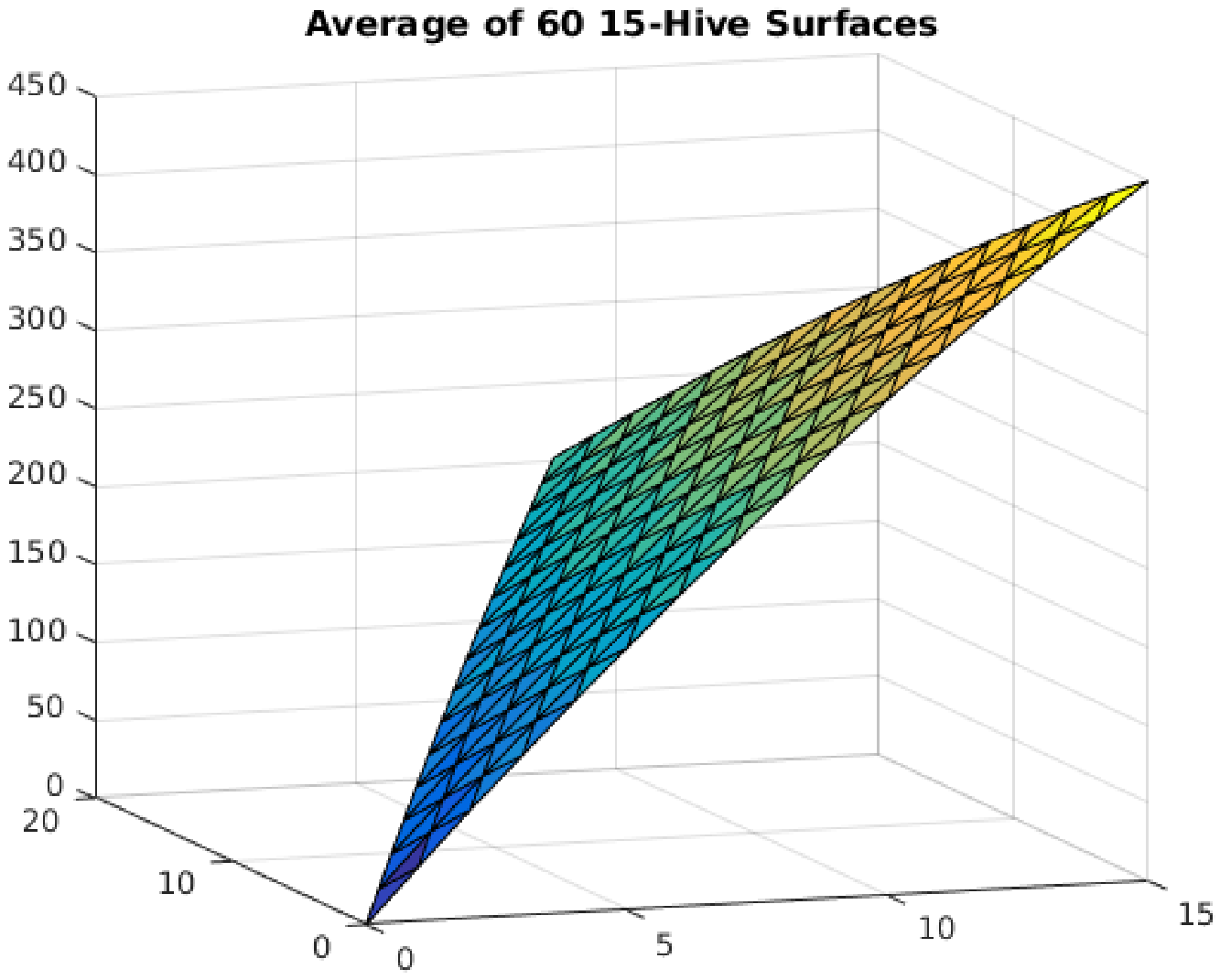}}
  \caption{Average hive surface on 60 identical pairs of normally distributed SPD(15) matrices on the left, and diagonally dominant SPD(15) on the right}
  \centerline{\includegraphics[scale = .6]{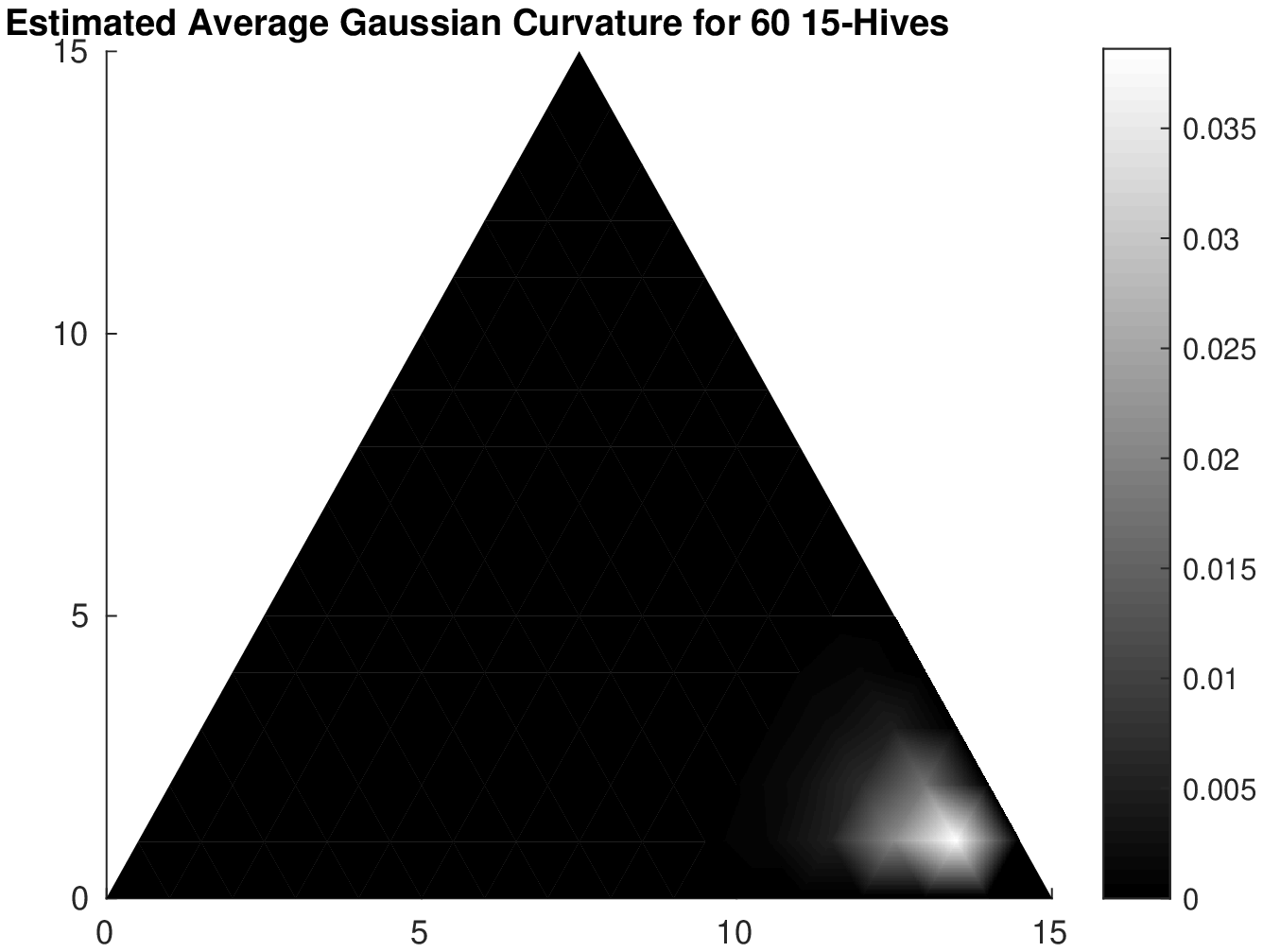} \includegraphics[scale = .6]{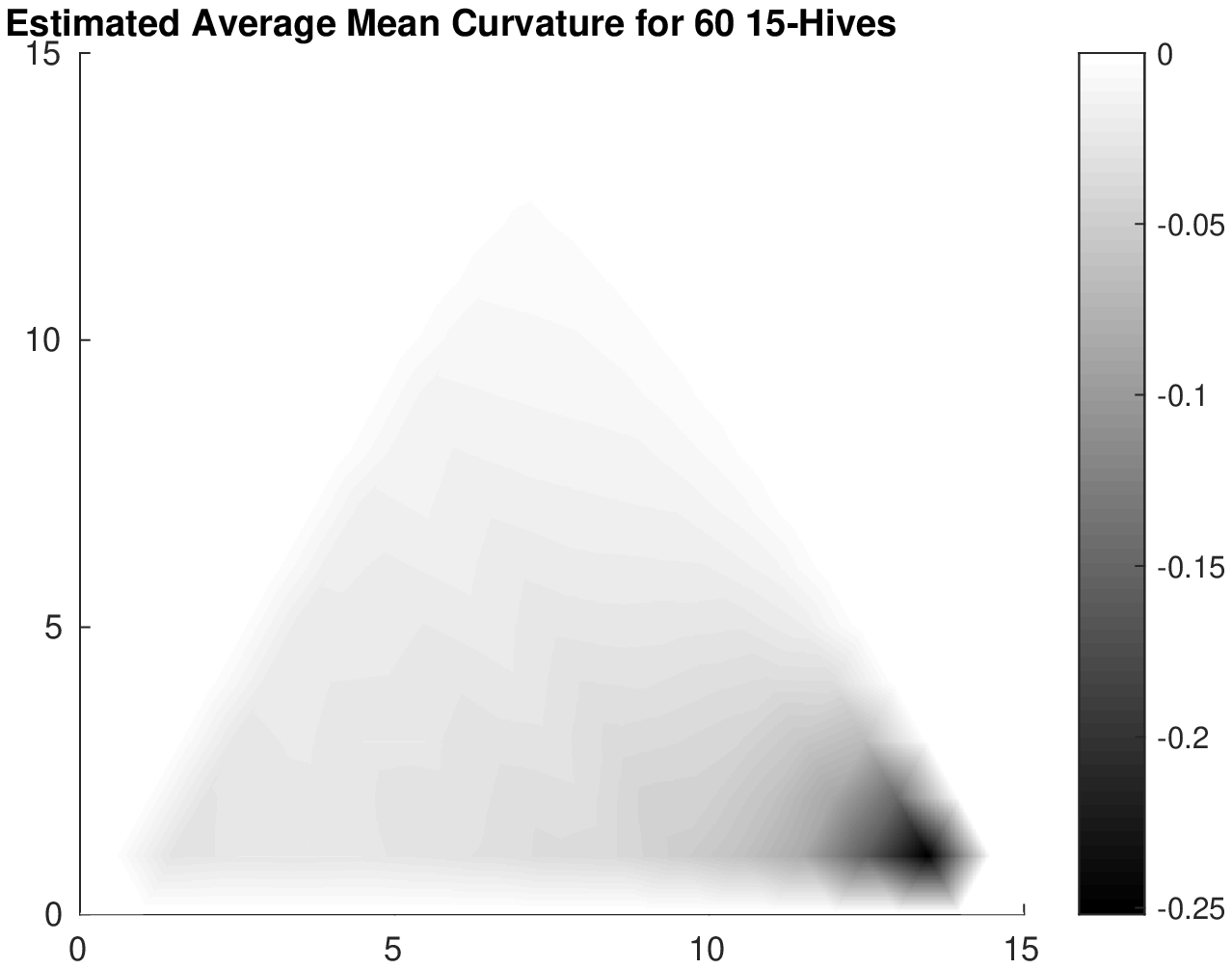}}
  \caption{Average gaussian and mean curvatures on 60 identical pairs of normally distributed SPD(15) matrices}
  \centerline{\includegraphics[scale = .6]{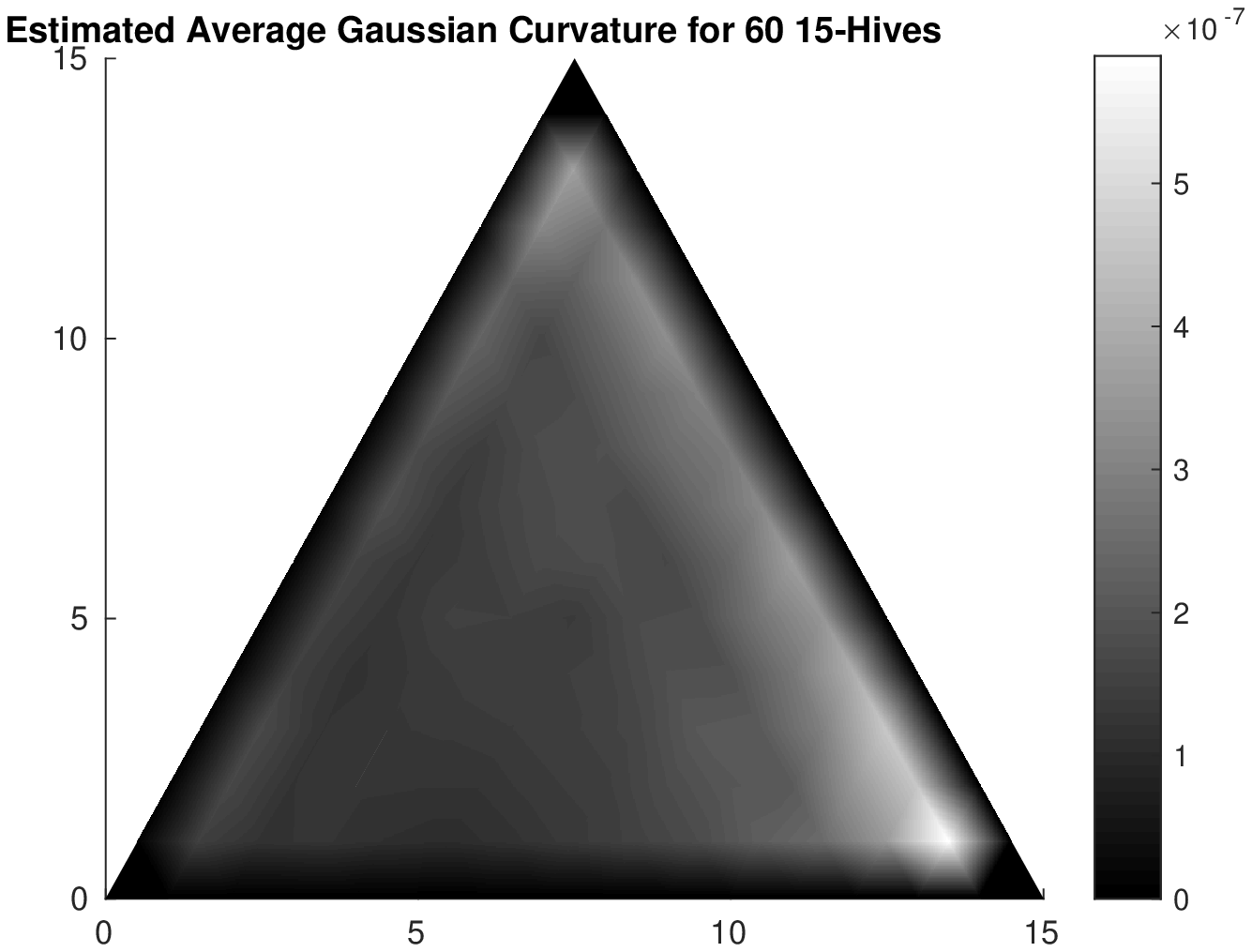} \includegraphics[scale = .6]{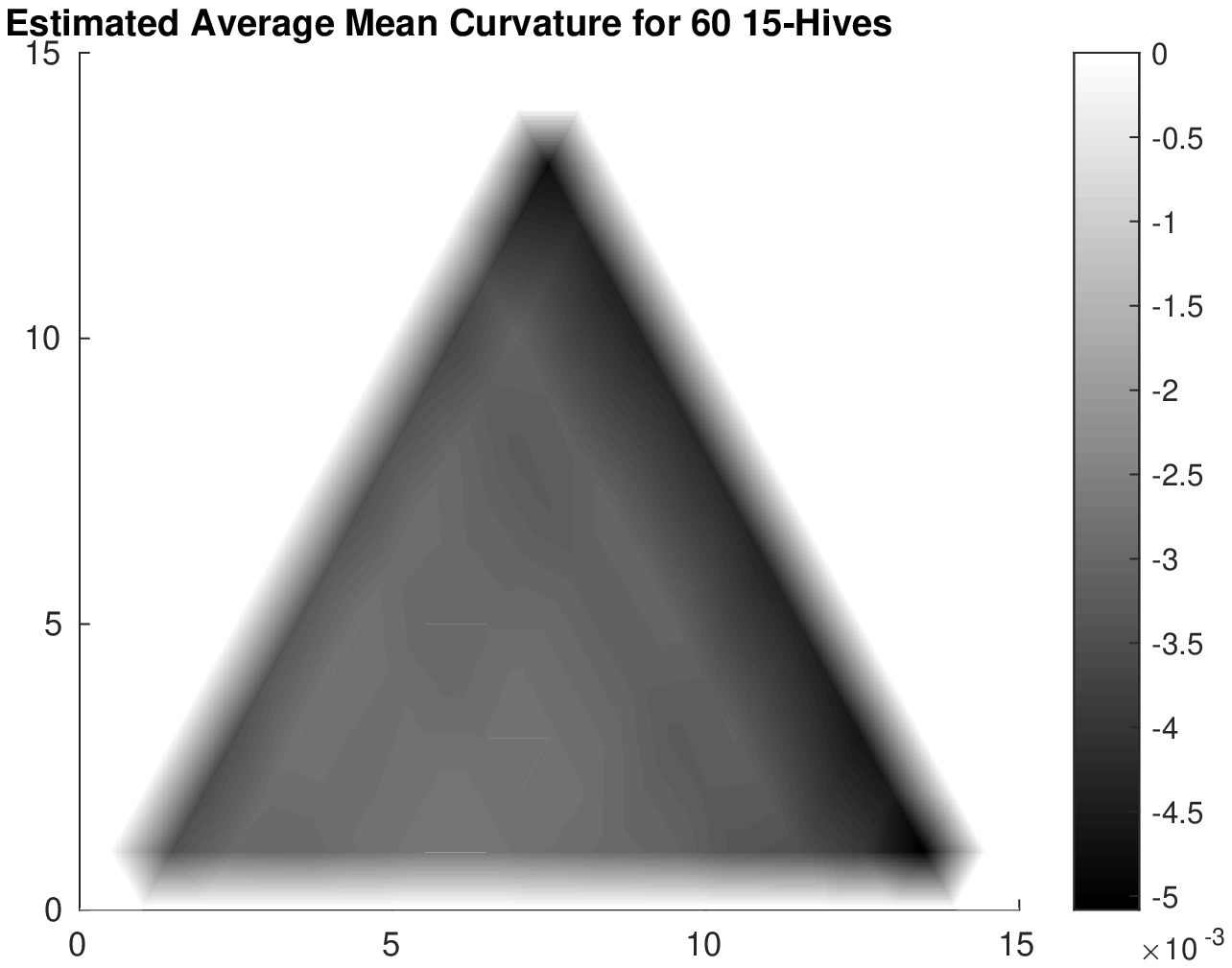}}
  \caption{Average gaussian and mean curvatures on 60 identical pairs of diagonally dominant SPD(15) matrices}
\end{figure}

\section{Summary of Hives from Hermitian Matrix Pairs}
The conjecture provided by Danilov and Koshevoy was the first generative scheme for finding interior hive coefficients in the sense of Knutson and Tao for pairs of Hermitian matrices. With the modification by Appleby and Whitehead which was accompanied by a proposed sketch of a proof, this optimization theoretically opened the door to a wide area of study of random continuum hives from various matrix ensembles, if the construction could be implemented as a numerical algorithm. We have analytically shown that the supporting arguments by AW are not generically true, which prevents such a general study from being undertaken using that construction. However, we provided the first implementation of their proposal by way of an analytic map onto a product space of Grassmannians that allowed for the AW construction to be tested, and we characterized the probability to which the optimization fails to find proper hives. Furthermore, we note several cases where we can generate hives with remarkable probability. In those ensembles, we have undertaken the first preliminary study of geometric observables of hive surfaces.

As noted, the usage of hives extends to the realm of computing Littlewood-Richardson coefficients for the integer valued cases. We now detail further results we have gained from studying combinatorial hives to that end.

\section{Rounded Estimation of LRCs using Hit-and-Run on the Hives} \label{sec:LRC}

It is possible to provide an estimate for the LRC associated with the 3-tuple of weight vectors $(\mu,\nu,\lambda)$ given the hive construction. In this section, we discuss our implementation of the algorithm described by Narayanan in ``Estimating Certain Non-Zero Littlewood-Richardson Coefficients'' \cite{Narayanan2014}. To our knowledge, this is the first implementation of a randomized approximation scheme for computing LRCs based on combinatorial hives. The algorithm relies on first forming the hive polytope $P_{\mu \nu}^{\lambda}$ constructed from the linear rhombus inequalities and the boundary equalities such that

\begin{equation}
  A x - b_{\mu \nu}^{\lambda} \preceq 0 \, ,
\end{equation}
where $A$ encodes the hyperplane structures and $b_{\mu \nu}^{\lambda}$ includes a vector of constants specifying the boundaries.

The polytope is then enlarged by a relaxation of the affine hyperplanes, providing a rescaling of the polytope $Q_{\mu \nu}^{\lambda}$ defined by some small fixed dilation as follows:
\begin{equation}
  A x - b_{\mu \nu}^{\lambda} \preceq 2 \, \,.
\end{equation}

An estimate for the continuum volume $\tilde{V}(Q_{\mu \nu}^{\lambda})$ of the enlarged polytope is then required. Once found, it is necessary to determine the number of integer lattice points of the interior original polytope which also fall inside the dilated polytope. Given this fractional estimate $f$, we can provide an approximation for the number of integer lattice points in the original hive as simply $V(HIVE(\mu,\nu,\lambda)) =  f \tilde{V}(Q_{\mu \nu}^{\lambda})$.

In our implementation, the enlarged polytope volume is computed by using Cousins and Vempala's ``A Practical Volume Algorithm,'' which includes empirical convergence tests, an adaptive annealing scheme and a new rounding algorithm built upon the framework of Lovasz and Vempala's nested convex body approach for volume estimation \cite{Cousins2016} \cite{Lovasz2006}. We take extra care in the preprocessing step of the polytope to avoid translations of the center of mass of the polytope to the origin, as this is not a lattice number preserving transform. However, we do perform the rounding preprocessing that violates the lattice volume. We discuss the implications of this in Section \ref{sec:preprocessing}.

We then perform a sampling of rational points in the larger polytope using an adaptive centering hit-and-run method with an implementation provided by Benham \cite{Benham2012}. We round the resulting points to the integer lattice, and test those points against the inequalities that define the original hive polytope, continuing to gather and test points until we achieve a convergence of the relative lattice volume estimate within some desired tolerance. Multiplying the outer continuum volume estimate by the fractional lattice volume estimate yields the desired integer lattice volume estimate for the original hive and thus the associated LRC.

\subsection{Results of the Rounded Algorithm} \label{sec:rounded}

\begin{figure}[H]
  \centerline{\includegraphics[scale = .5]{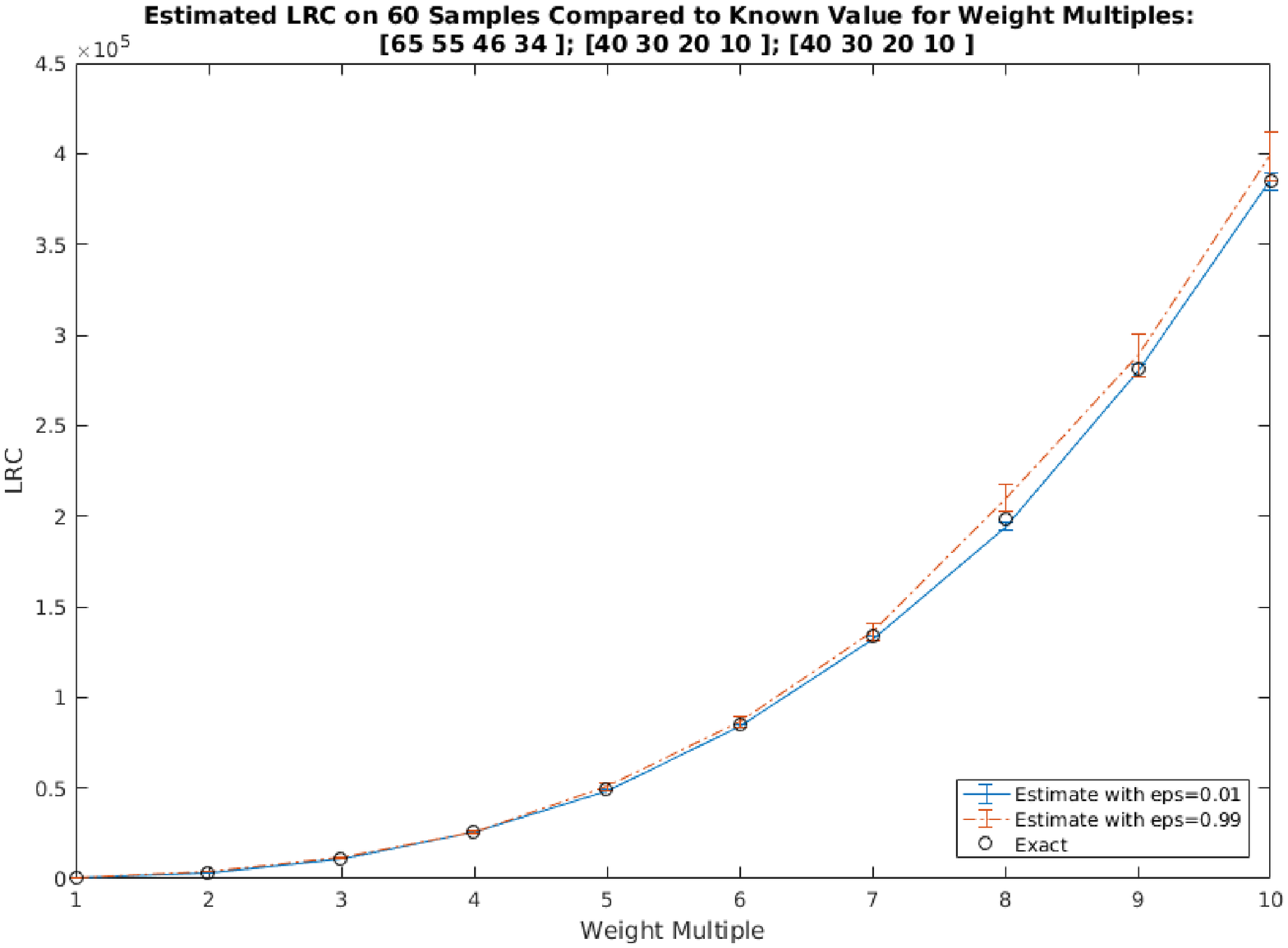}}
  \centerline{\includegraphics[scale = .5]{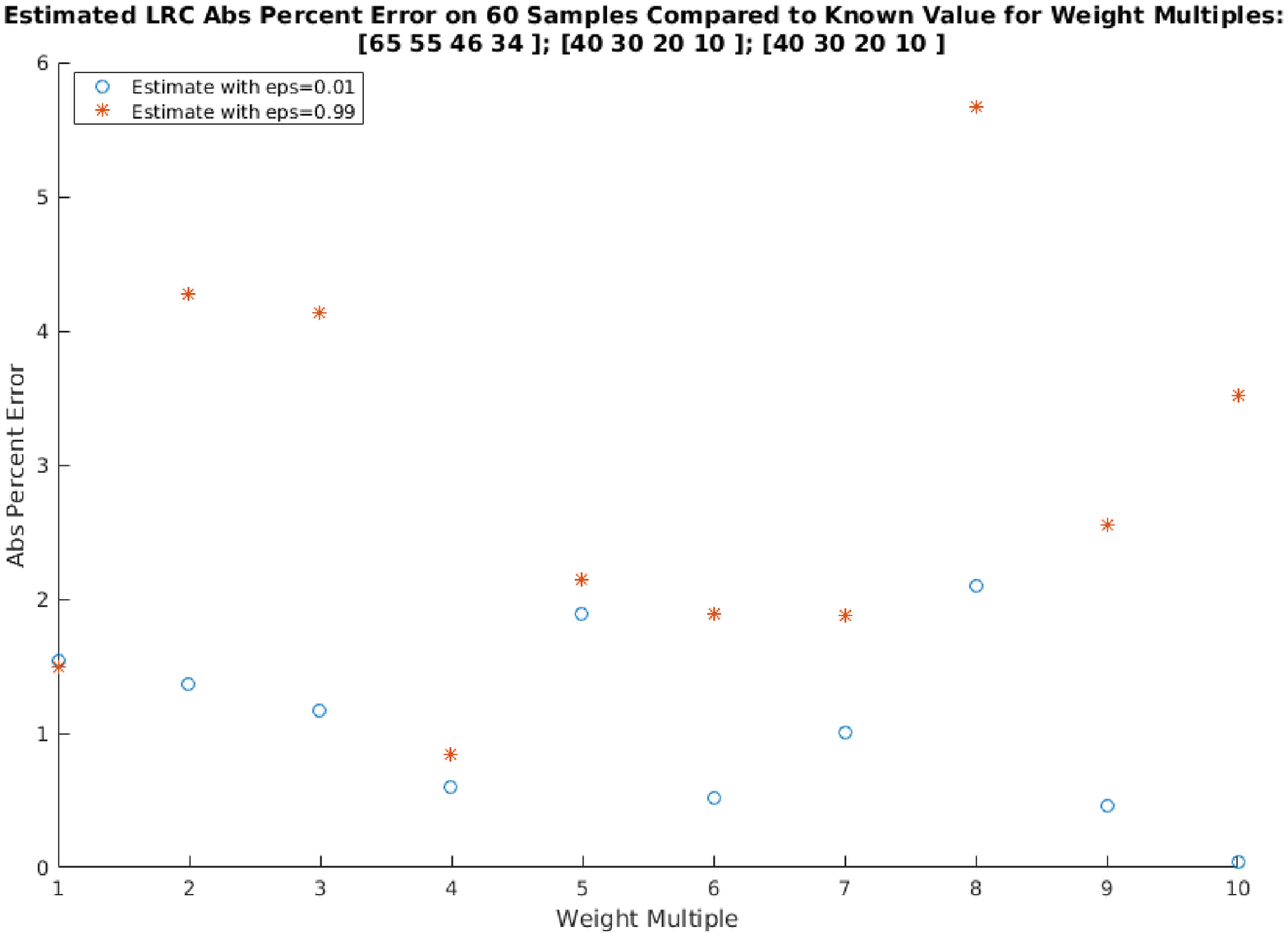}}
  \caption{Accuracy of rounded LRC estimator for weight vector multiples, 95\% confidence intervals and absolute percent error} \label{fig:CLRCCI}
\end{figure}

We test the accuracy of this algorithm against known LRC values with respect to a user defined relative error parameter. The known LRC values were computed using Anders Buch's LRC C-based calculator which uses a computationally expensive exhaustive algorithm \cite{Buch1999}. We pick a known 4 dimensional tuple, and look at scalar multiples of these weights that correspond to dilations of the hive polytope as a way of systematically testing the robustness of the algorithm for finding coefficients from progressively larger convex bodies. We see in Fig. \ref{fig:CLRCCI} that the absolute accuracy of the algorithm decreases with hive volume, but the lower relative error parameter is able to compensate at the cost of computational time. Regardless of the error parameter, the algorithm has un-patterned fluctuations in its relative percent error with a wide deviation and no clear trend as a function of volume. However, we were able to achieve less than 5\% absolute relative error out to LRC's with hive volumes on the order of $10^5$ in a under $5$ seconds of computation time on a variety of modern architectures, where the exact Buch implementation would take on the order of minutes.

\section{Coordinate Hit-and-Run on the Hive Lattice Itself}

We created a second estimation algorithm for the hive volumes, where instead of working in the continuum and rounding to the integer lattice, we attempt to estimate the lattice volumes directly using a lattice coordinate hit-and-run (CHAR). We use a similar volume ratio approach, but work directly on the lattice and utilize smart resampling in the style of Ge, Ma, and Zhang's ``A Fast and Practical Method to Estimate Volumes of Convex Polytopes'' \cite{Ge2014}. Our algorithm takes the following form.

\begin{figure}[H]
  \centerline{\includegraphics[scale = .5]{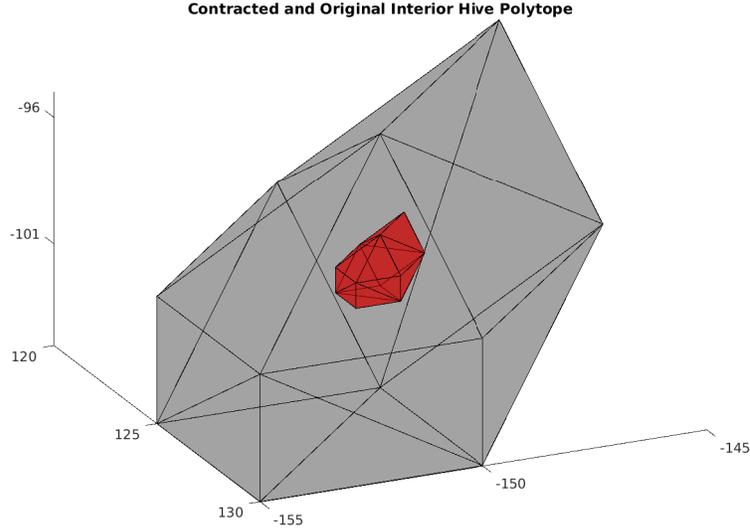}}
  \caption{Original and contracted hive polytope for a 4D weight vector tuple}\label{fig:hivePolytope}
\end{figure}

First, we formulate the maximum LP from the boundary data as it is guaranteed to give an integer point in the hive by Knutson and Tao. Assuming one can be found, this lattice point will be the start of a CHAR algorithm that is designed to approximate the full lattice volume of the hive. We avoid performing any rounding of the hive polytope, as this would not preserve lattice volume. Since the rhombus inequalities are constraints on 4 interior hive coordinates for each possible rhombus, the hyperplane matrix $A$ has at most 4 non-zero entries per row with values $\pm 1$. This indicates that the polytopes themselves are not markedly dissimilar to hypercubes in terms of their aspect ratios, and we do not expect extreme channels of small volume to appear that would cause conductance issues in our sampling algorithm. An illustration of this can be seen in Fig. \ref{fig:hivePolytope}, where we show the polytope corresponding to 4D weight vectors that yield the 3D interior hive volume from the 3 free bulk coordinates. The structure takes the form of a multi-faceted rounded rectangle that does not have an extreme aspect ratio or have narrow regions separating large volumes, etc.

With our initial integer hive found, we then progressively constrain the hyperplanes uniformly in the LP by solving a series of LPs

\begin{equation}
    A x - b_{\mu \nu}^{\lambda} \preceq \xi
\end{equation}

for $\xi \in [0,\ldots,\xi^*] \subset \mathbb{Z}^-$, shrinking the hive volume until no solution can be found at $\xi^*$, which indicates that there is no longer an integer point in the interior of the hive volume. We take as a shrunken polytope the hive defined by
\begin{equation}
    A x - b_{\mu \nu}^{\lambda} \preceq \min{(0,\tilde{\xi}\equiv \xi^* + 2)}
\end{equation}
in order to avoid excessively small volumes that may induce lattice points that are not accessible by the coordinate alignment of our hit-and-run and give a false lattice volume. This contracted polytope is also illustrated in Fig \ref{fig:hivePolytope}.

At this point, the algorithm takes a branching path: if the starting `contracted' polytope is the original, the number of lattice points should be small enough that a stochastic enumerative sampler can be utilized to calculate the hive volume. That is, we can perform a CHAR on the original hive saving only unique lattice points until the number of such unique points appears to become stationary with respect to some computational heuristic. The number of discovered lattice points directly yields the LRC for the weight tuple.

If the hive volume is large, multiple contractions will be necessary. In this branch, we need to estimate the lattice size, as direct counting methods are clearly unfeasible (especially in high dimensions, where the curse of dimensionality will prevent us from reliably sampling from a subset of the target space). Given $\tilde{\xi}$, we now have a sequence of nested polytopes with decreasing volume from our desired original hive such that each contracted polytope is wholly contained in the next larger, and the geometry is preserved (in contrast to the methods in \cite{Ge2014} that rely on nested hyperballs within the desired convex body which require delicate arguments and rounding procedures to ensure that in high dimensions the bulk of the desired polytope volume is not excised by even the largest inscribing sphere). Following the volume ratio arguments provided by \cite{Narayanan2014}, we consider the sequence of $n$ polytopes $P_{\mu \nu}^{\lambda} (\xi_n \in [\tilde{\xi},\tilde{\xi}+2,\ldots,0])$ such that $P_{\mu \nu}^{\lambda} (\xi_0) = P_{\mu \nu}^{\lambda}(\tilde{\xi})$ and $P_{\mu \nu}^{\lambda} (\xi_n) = P_{\mu \nu}^{\lambda}(0) = P_{\mu \nu}^{\lambda}$. This sequence admits volume ratios between successive contractions which scale by at most a constant. If $\tilde{\xi}$ is not even, we include one additional ratio down to the minimum that will clearly have a volume ratio even closer to 1 than the others. We assume $\tilde{\xi}$ is even for this discussion.

By the telescoping volume argument,

\begin{equation}
  vol(P_{\mu \nu}^{\lambda}) = vol(P_{\mu \nu}^{\lambda}(\tilde{\xi})) \prod_{n = 0}^{n-1}{\frac{vol(P_{\mu \nu}^{\lambda}(\xi_{n+1}))}{vol(P_{\mu \nu}^{\lambda}(\xi_n))}} \, .
\end{equation}

We push this to the $d$-dimensional lattice volumes contained within each convex contraction, letting our hive sequence be defined by $HIVE(\mu,\nu,\lambda) (\xi_n) = P_{\mu \nu}^{\lambda} (\xi_n) \cap \mathbb{Z}^d $.

The innermost hive is estimated directly by a unique accumulating CHAR just as in the minimal case. We then estimate the successive lattice volume ratios by first performing a fixed number of steps in a CHAR on the largest (original) hive. We check to see how many of these points fall inside the next inner contracted polytope. The ratio of the total number of samples to the number of samples that fall inside the next contracted interior is proportional to our desired volume ratio. To save computational time, we note that the original random walk on the largest space, if enough data was sampled for the distribution to be nearly uniform, will induce a nearly uniform distribution of points in each of the contracted hives as well. As a result, we can then save the points that tested positive to lie in the interior of the contracted hive, and resample from the smaller hive only enough points such that the relative error between each ratio remains constant. This outer-to-inner sampling procedure was demonstrated by \cite{Ge2014} to save over 70\% time consumption for a wide range of polytopes. Repeating this process until all of the lattice volume ratios are known will yield the required telescoping product constituents. For the inner starting points, we choose from the uniform distribution over the previously sampled points known to be in the contracted interior, and proceed with a new CHAR.

Our implementation of the lattice CHAR works directly with the rhombus constraints. The maximal hive produced by the LP serves as a starting hive for the random walk. We treat each interior bulk hive index as a coordinate direction in the polytope space. We systematically check each index for `flex': that is, until we find a coordinate that has the ability to take on a different integer value while still satisfying the rhombus inequalities, we proceed to examine the 12 rhombuses associated with each coordinate (4 vertical, 4 left facing, 4 right facing) and thereby find the integer subset over which the inequalities would still hold. This is an exact constant time process unlike Cousins' algorithm that, although general to handle arbitrary convex bodies, requires a bisection subdivision search to approximate the intersection of the hit-and-run axis with the boundary of the polytope. If each coordinate has an empty subset, the hive is `tight' and may indicate that there is only 1 integer solution for the provided weight vectors--the algorithm then returns. The first instance of flex becomes the first coordinate choice, wherein a new value is chosen from the uniform distribution on the available integer line segment. New coordinates are chosen uniformly randomly from the total bulk index set, their allowed 1D integer spaces are computed, and new values are selected. The random walk proceeds in this way for a fix number of steps, ensuring that we always remain in the HIVE space by construction (therefore no rejection sampling is needed), and we walk directly on the integer lattice without having to round our walk.

\subsection{Results of the Lattice Algorithm}
We note that this algorithm relies on a few key assumptions. Foremost, we require that our inner hive estimation is good. A variety of convergence heuristics can be employed to determine when the inner unique walk has settled; however, if the contracted volume has lattice points that are not aligned with our coordinatization of the hive polytope, then this process will always under-represent the hive volume. For example, each individual coordinate may appear tight with respect to the starting solution from the LP. However, global changes in the bulk coordinates (or the simultaneous change in simply more than one coordinate at a time) may have enabled a new integer point to be found. This is not accessible from our CHAR due to our basis choice. Although it is a natural basis for the hive lattice, there is no guarantee that the polytope is coordinate aligned.

Future work in characterizing unimodular (lattice preserving) transforms of the hive polytope in order to better align the lattice with the hive coordinates is underway.

\begin{figure}[H]
  \centerline{\includegraphics[scale = .5]{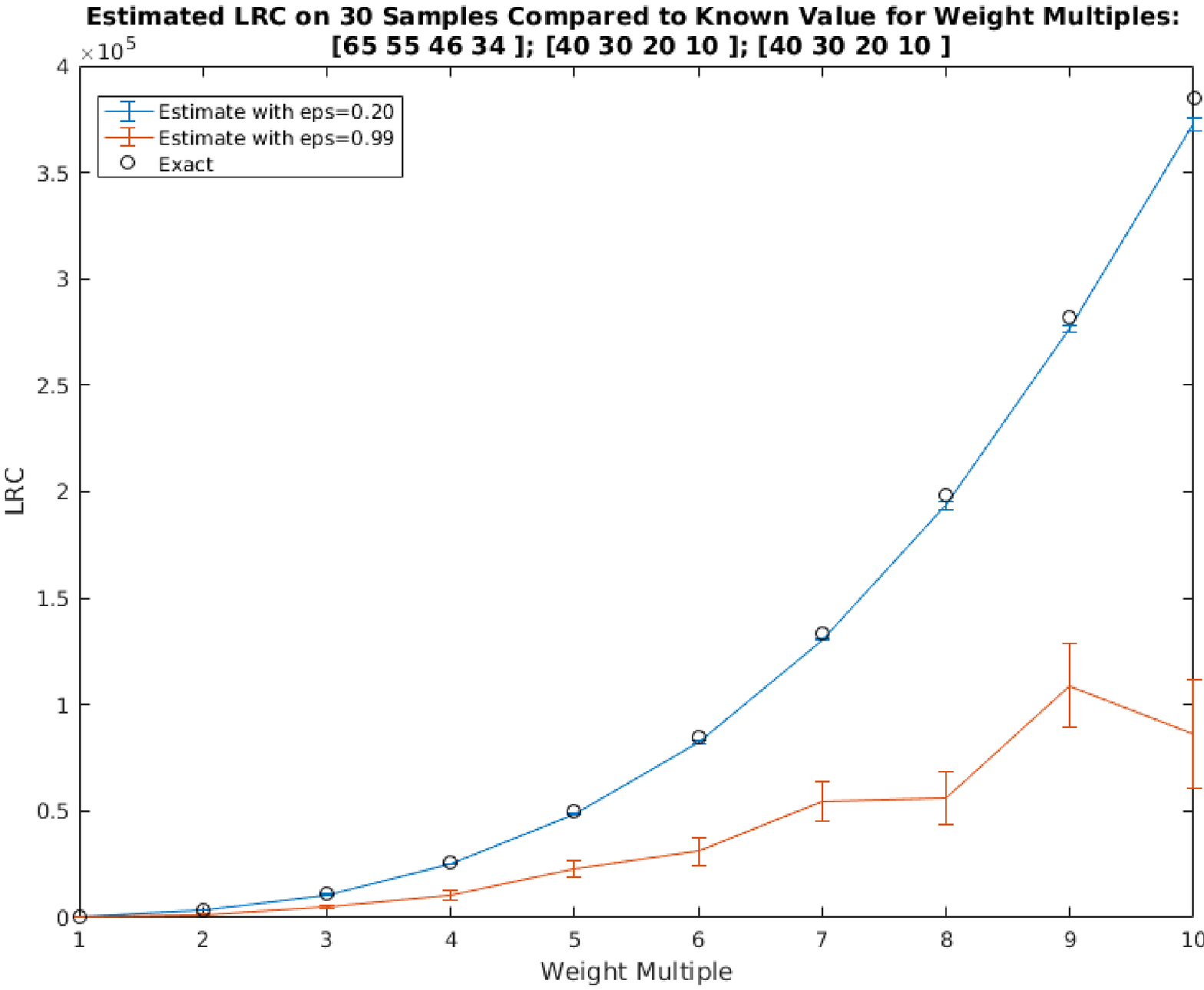}}
  \centerline{\includegraphics[scale = .5]{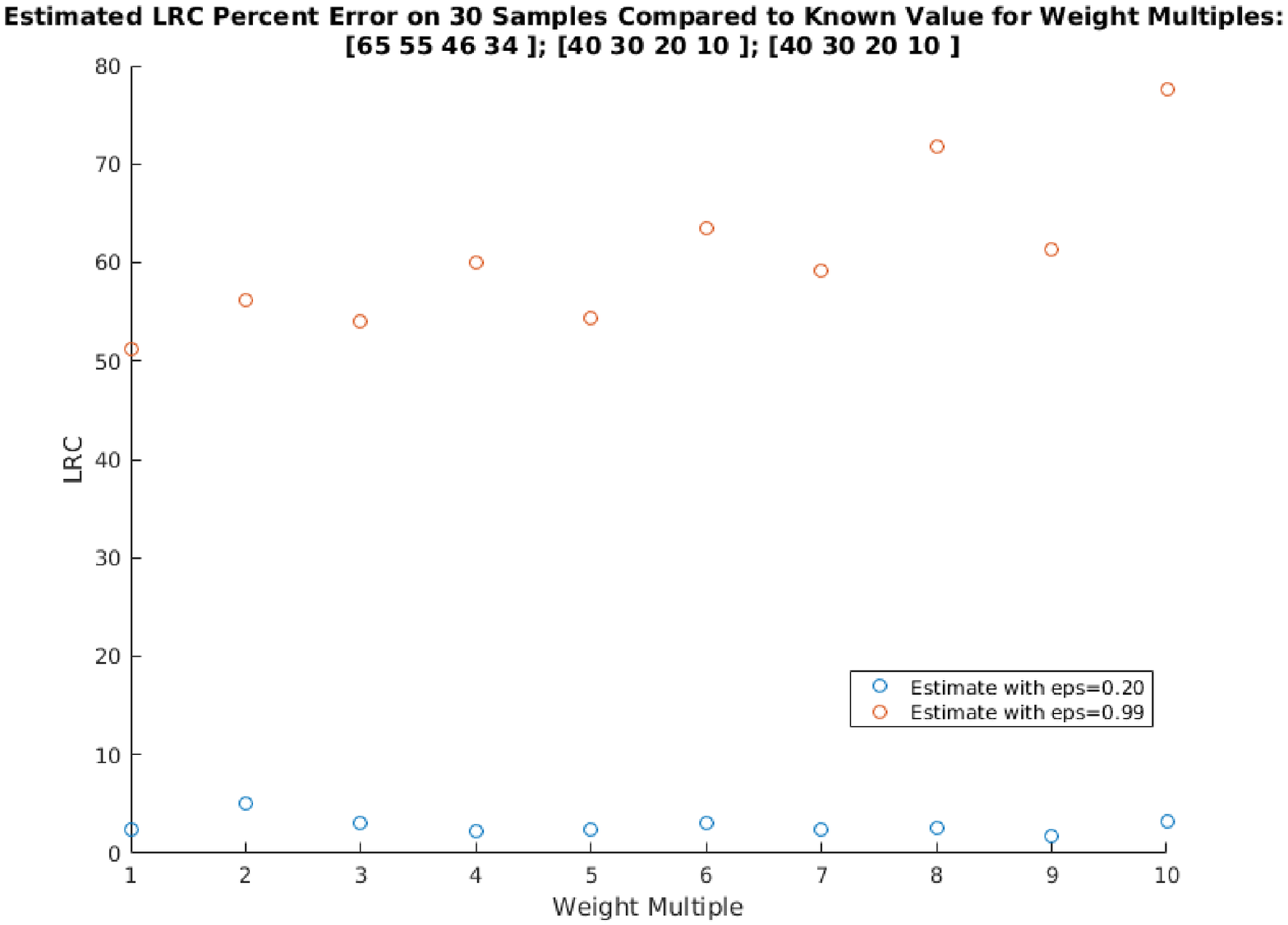}}
  \caption{Accuracy of lattice LRC estimator for weight vector multiples, 95\% confidence intervals and absolute percent error} \label{fig:LRCCI}
\end{figure}

Fig. \ref{fig:LRCCI} shows the accuracy of this algorithm in the same comparison test used above with the Cousins' volume estimation. We see a large degradation of accuracy with the large error parameter. This is expected, as the error parameter conditions the initial convergence for the unique point accumulator in the minimal hive volume search, and if this is inaccurate, all subsequent products against this ratio will will be inaccurate. Lower relative error parameters product results competitive with the rounded algorithm described in Sec. \ref{sec:rounded}.

\begin{figure}[H]
  \centerline{\includegraphics[scale = .55]{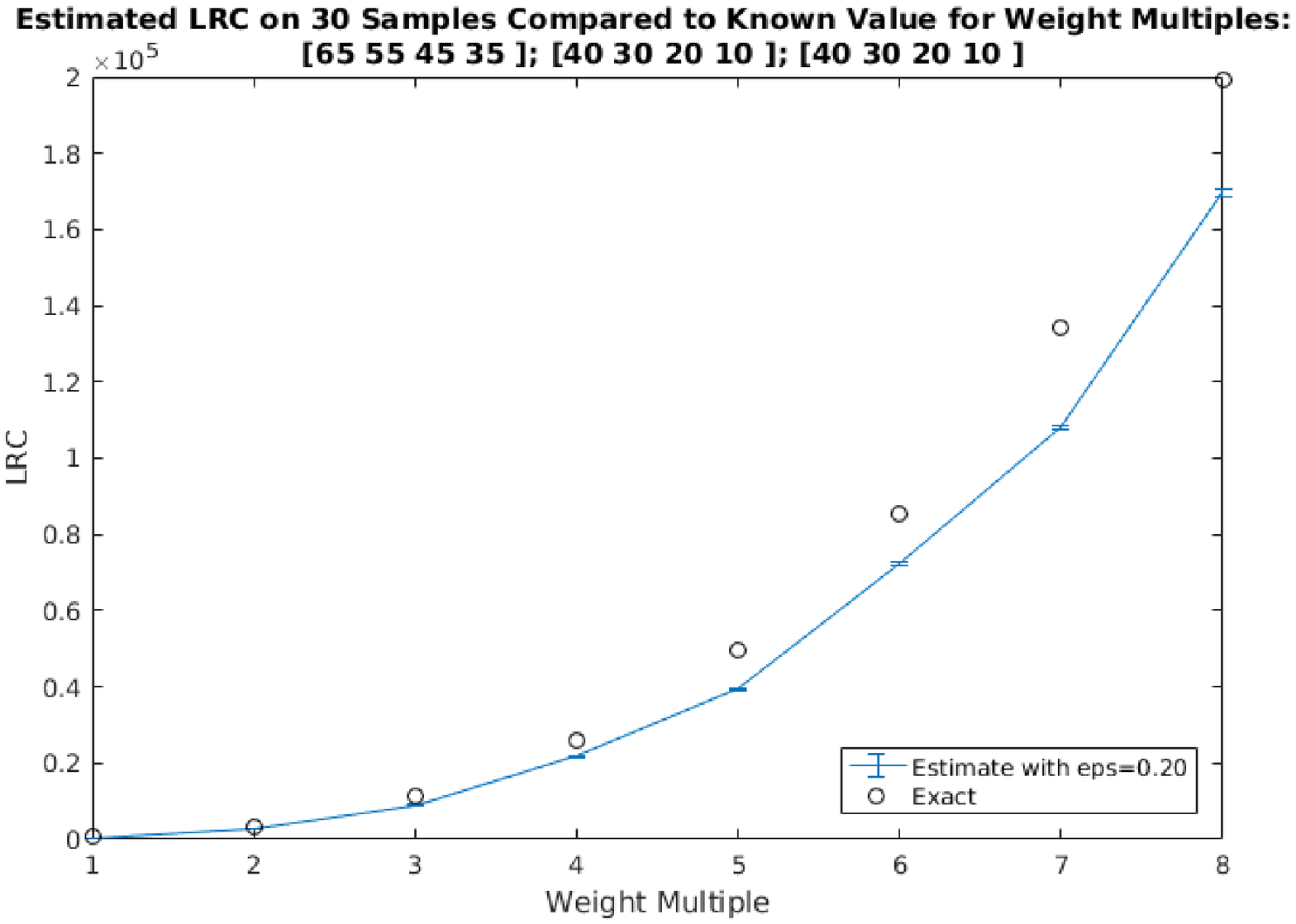}}
  \centerline{\includegraphics[scale = .55]{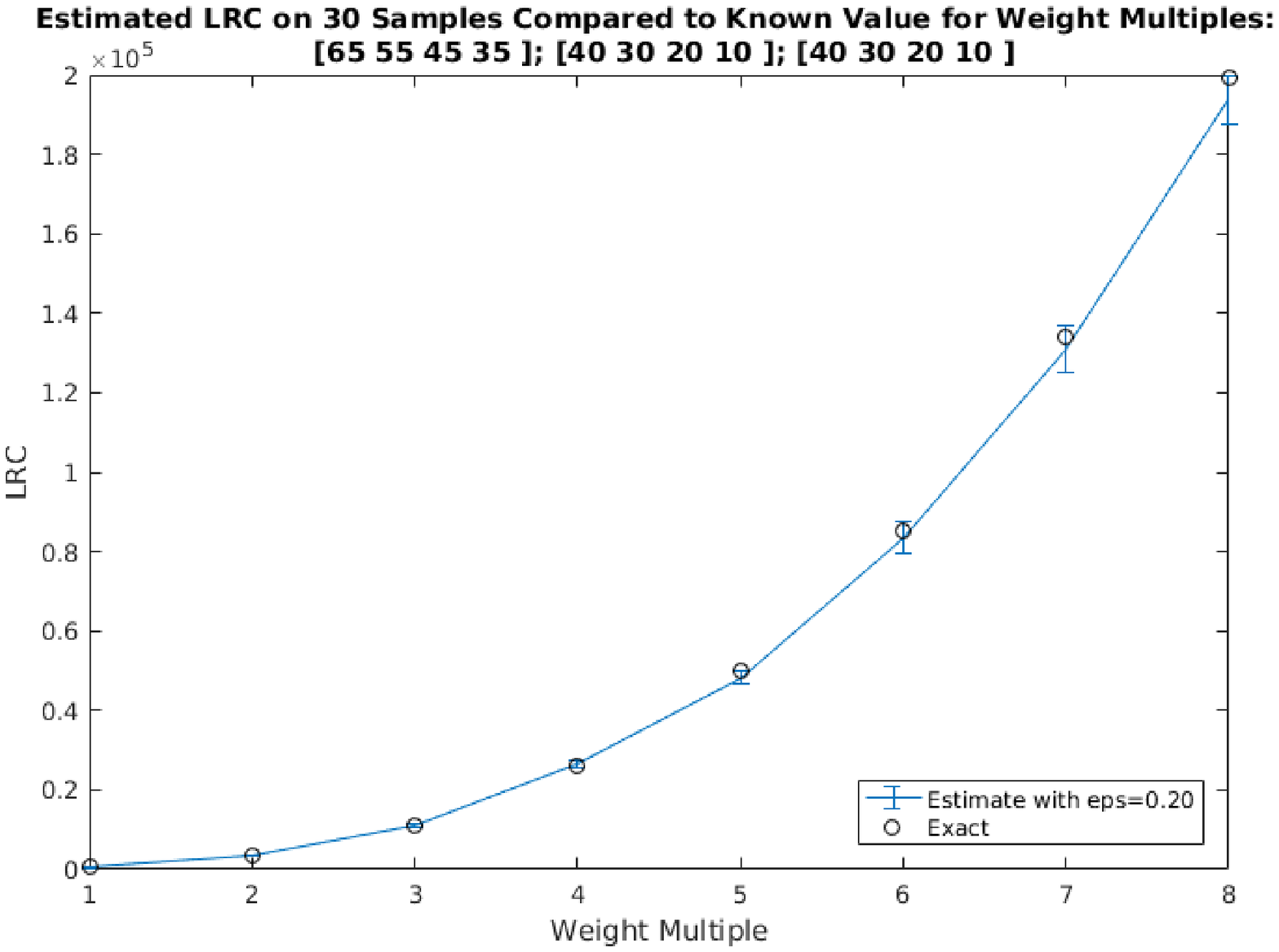}}
  \caption{Accuracy of LRC estimators for weight vector multiples with poor lattice alignment with the rounded estimator on the bottom} \label{fig:LRCMultiplesAccBad}
\end{figure}

The above results were hive volume scalings wherein the algorithm was able to perform a connected walk over the full hive space. A unique accumulating run over the original full lattice would yield the exact LRC of 505 for the base weight vectors. However, the tuple $([40,30,20,10],[40,30,20,10],[65,55,45,35])$, which is a near perturbation of the other tuple and has an LRC of 506, appears to contain a lattice point that is not aligned with the basis of our CHAR, and as a result, any long-running accumulating sampler only finds 505 unique configurations. This has dramatic effects on the quality of the estimations in general once we consider multiples of the weight vectors, and one can see the difference in the fractional errors illustrated in Fig. \ref{fig:LRCMultiplesAccBad} for the same low error parameter as in the above sampling compared with other volume estimation algorithm over the same set. Now, even for the previously acceptable error parameter for roughly the same sized polytope volumes, we have poor accuracy. As mentioned, unimodular transforms may aid in this type of systematic inaccuracy, but this may be a characteristic failure of fixed axis CHAR algorithms on embedded discrete sets.

Although there are no proven bounds on the mixing time of CHAR algorithms in the continuum, in practice they have been shown to have the same properties as traditional hit-and-run, and with fewer computational steps \cite{Emiris2014}. On lattice walks, however, the CHAR is strongly susceptible to the issues presented here unless appropriate care is taken based on foreknowledge of the lattice volume within the convex body one is trying to sample.

\begin{figure}[H]
  \centerline{\includegraphics[scale = .55]{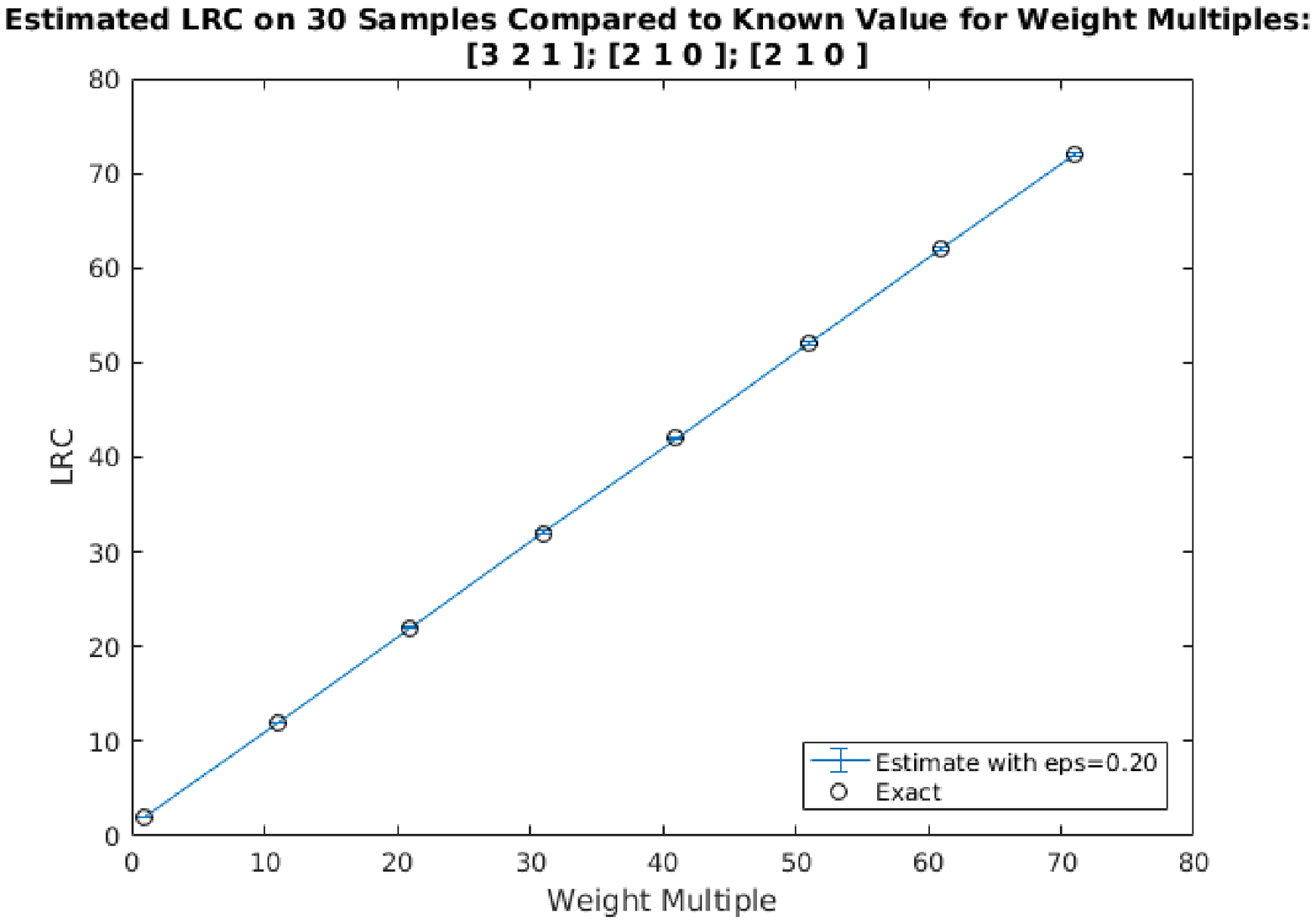}}
  \centerline{\includegraphics[scale = .55]{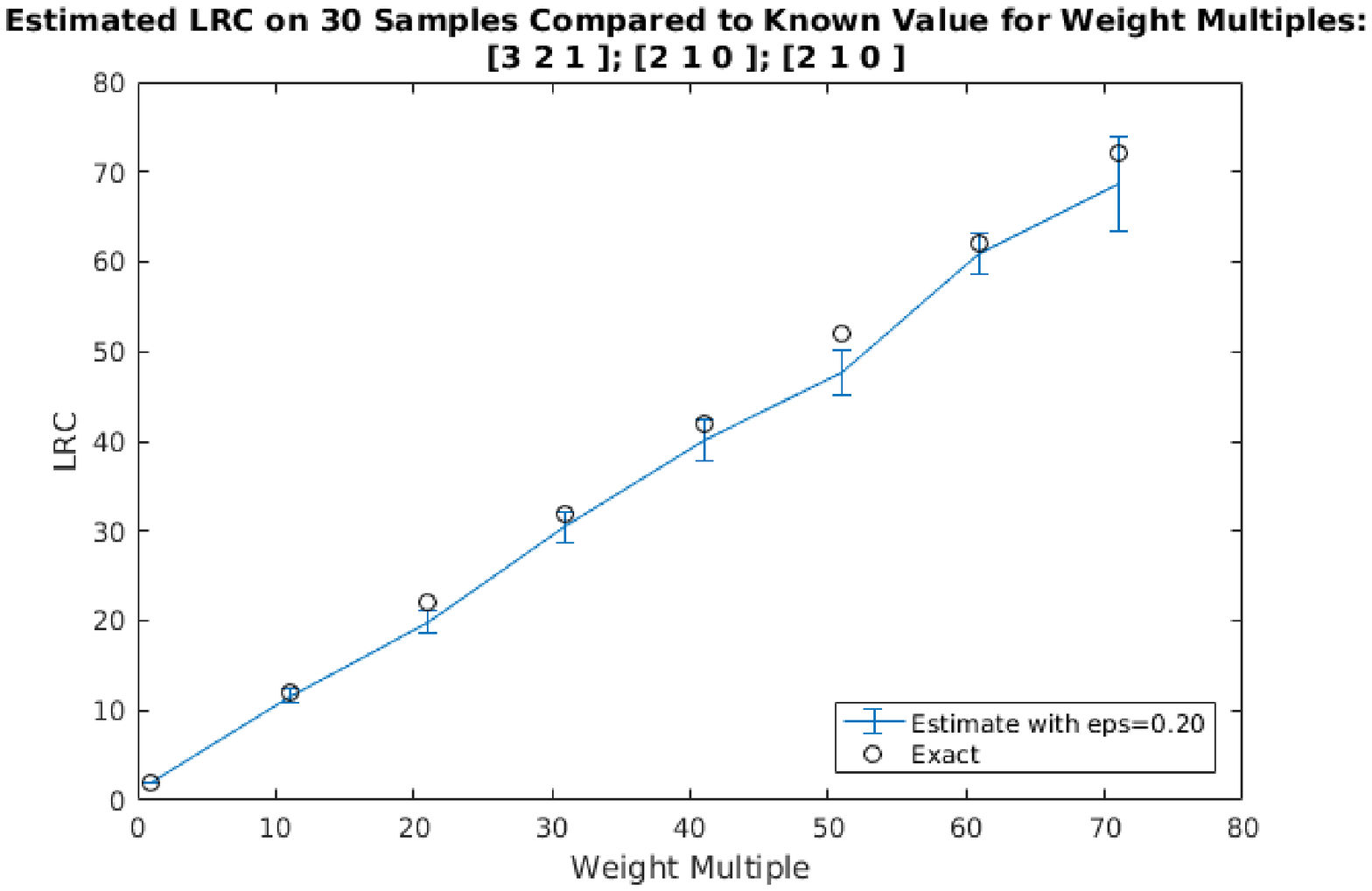}}
  \caption{Accuracy of LRC estimators for weight vector multiples with small volume with the rounded estimator on the bottom} \label{fig:LRCMultiplesAccSmall}
\end{figure}

In the limit where the number of lattice points is assumed to be large, the rounding pre-processing of the polytope will not affect the estimation on an order more significant than the errors from the volume estimation itself, despite the fact that the transforms are not lattice preserving\label{sec:preprocessing}. However, for small lattice volumes, this affect can be comparable, and the direct lattice walk that preserves the structure outperforms in accuracy as shown in Fig \ref{fig:LRCMultiplesAccSmall}.

\section{Future Work}
We have presented an exploration of a few analytical and numerical aspects of combinatorial hives and their applications. Future directions for this research would involve finding an appropriate modification of the Appleby-Whitehead construction to produce real valued hives in all cases, or to that end, re-examining the original conjecture by Danilov and Koshevoy and probing their proposal numerically using similar techniques to our Grassmannian decomposition and testing for its validity. With respect to Littlewood-Richardson calculations, finding and implementing unimodular rounding schemes would be greatly beneficial for both algorithms, producing more accurate results with respect to lattice counts. A thorough study of the algorithms for much higher dimensional vector spaces should also be undertaken to better characterize their efficiency.

\section{Acknowledgements}
We would like to thank Hariharan Narayanan for the introduction to this construction and for many helpful conversations as an advisor and mentor. Additional thanks to Glenn Appleby in corresponding about \cite{Appleby2014} and clarifying the history of hives from Hermitian matrix pairs with respect to reference \cite{Danilov2003}. We would also like to thank Github for a student developer grant and the hosting of private repositories. The algorithms discussed in this paper will be made freely available for public use and modification upon publication. This work was supported in part by the University of Washington.

\appendix*
 
\section{Euclidean Computations for the Grassmannian Map}

To proceed with a trust-region or gradient descent method, we need to compute the Euclidean gradients and Hessians of our cost function:

\begin{eqnarray}
  H_{ijk} &=& \min_{\substack{\tilde{A}=\text{col}(\vec{\alpha}_i) \,|\,\text{span}(\vec{\alpha}_i) \perp U \\ B=\text{col}(\vec{\beta}_i) \,|\,U=\text{span}(\vec{\beta}_i) \subseteq V}}
  -(f(B,\tilde{A},M)+g(B,N)) \nonumber \\
  f(B,\tilde{A},M) &\equiv& (tr(\rv{B}{\tilde{A}}(\rv{B}{\tilde{A}}^T\rv{B}{\tilde{A}})^{-1}\rv{B}{\tilde{A}}^T M \rv{B}{\tilde{A}}(\rv{B}{\tilde{A}}^T\rv{B}{\tilde{A}})^{-1}\rv{B}{\tilde{A}}^T)  \\
  g(B,N)&\equiv& tr(B(B^TB)^{-1}B^T N B(B^TB)^{-1}B^T)) \nonumber
\end{eqnarray}
for functions $f$ and $g$ defined as the respective portions of the cost function.

\subsection{Gradients}

\subsubsection{$\nabla_{\cdot} f$}

  We compute the first order Frechet derivative of $f$ with respect to $B$ to begin, using a small perturbation $d$.
  Recall that the Frechet derivative of a function $f\,:\,S\rightarrow \mathbb{R}$ at a point $s$ in a normed linear space $S$ is defined as the unique linear operator $\mathscr{D}$ tangent to $f$ at $s$ that satisfies
  \begin{equation}
    f(s+d) = f(s)+\mathscr{D}(d) + \mathcal{O}^+ \,.
  \end{equation}

  On the space of real matrices equipped with the Frobenius norm, the gradient should satisfy the relation
  \begin{equation}
    \mathscr{D}(d) = \langle \nabla_s f,d \rangle = tr((\nabla_s f)^T d) \,.
    \end{equation}

  We use the notation $[left]$ to denote duplication of `the same block of terms as expanded on the left,' for conciseness.

  \begin{eqnarray}
    f(B+d,\tilde{A},M) &=& tr(\rv{B+d}{\tilde{A}} (\cv{B^T+d^T}{\tilde{A}^T}\rv{B+d}{\tilde{A}})^{-1} \cv{B^T+d^T}{\tilde{A}^T} M [left]) \nonumber  \\
    &=& tr(\rv{B+d}{\tilde{A}} (\tbt{B^TB}{B^T\tilde{A}}{\tilde{A}^TB}{\tilde{A}^T\tilde{A}} + \tbt{B^Td+d^TB}{d^T\tilde{A}}{\tilde{A}^Td}{0})^{-1} \cv{B^T+d^T}{\tilde{A}^T} M [left]) + \mathcal{O}^+
  \end{eqnarray}

  Let
  \begin{equation}
    X \equiv \tbt{(X_1)_{k\times k}}{(X_2)_{k\times i}}{(X_3)_{i\times k}}{(X_4)_{i\times i}}\equiv (\rv{B}{\tilde{A}}^T\rv{B}{\tilde{A}})^{-1}
  \end{equation}
  be the block matrix representing the inverse of the symmetric matrix $\rv{B}{\tilde{A}}^T\rv{B}{\tilde{A}}$.
  Given that inverses preserve symmetry, we note the following properties:
  \begin{eqnarray}
    X^T &=& X \nonumber \\
    X_1^T &=& X_1 \nonumber\\
    X_4^T &=& X_4 \\
    X_2^T &=& X_3\nonumber
  \end{eqnarray}

  Using this notation, we can expand the inverse assuming our perturbation $d$ is infinitesimal:

  \begin{equation}
    f(B+d,\tilde{A},M) = tr(\rv{B+d}{\tilde{A}} (X - X \tbt{B^Td+d^TB}{d^T\tilde{A}}{\tilde{A}^Td}{0} X) \cv{B^T+d^T}{\tilde{A}^T} M [left]) + \mathcal{O}^+ \,.
  \end{equation}

  Expanding all terms left of $M$ up to first order in $d$, we result in two terms built out of sums of $n\times n$ matrices that we have grouped according to their sign:

  \begin{eqnarray}
    f(B+d,\tilde{A},M) &=& tr([ \{BX_1B^T + d X_1 B^T + \tilde{A}X_3B^T + BX_1d^T+\tilde{A}X_3d^T + BX_2\tilde{A}^T + \tilde{A}X_4 \tilde{A}^T \} \nonumber \\
    &-& \{BX_1B^TdX_1B^T + BX_1d^TBX_1B^T+BX_1d^T\tilde{A}X_3B^T + \tilde{A}X_3B^TdX_1B^T+\tilde{A}X_3d^TBX_1B^T+\tilde{A}X_3d^T\tilde{A}X_3B^T \nonumber \\
    &+& BX_2\tilde{A}^TdX_1B^T + \tilde{A}X_4\tilde{A}^TdX_1B^T + BX_1B^TdX_2\tilde{A}^T+BX_1d^TBX_2\tilde{A}^T + BX_1d^T\tilde{A}X_4\tilde{A}^T \nonumber \\
    &+& \tilde{A}X_3B^TdX_2\tilde{A}^T + \tilde{A}X_3d^TBX_2\tilde{A}^T + \tilde{A}X_3d^T\tilde{A}X_4\tilde{A}^T + BX_2\tilde{A}^TdX_2\tilde{A}^T + \tilde{A}X_4\tilde{A}^TdX_2\tilde{A}^T \} ] M [left]) \nonumber \\
    &+& \mathcal{O}^+
  \end{eqnarray}

  We now introduce useful shorthands to illustrate some of the symmetry in this equation.
  Let the second term in curly brackets be denoted $\{*\}$. Similarly, let
  \begin{equation}
    \{**\} \equiv \{dX_1B^T + BX_1d^T + \tilde{A}X_3d^T = dX_2\tilde{A}^T - \{*\} \} \, ,
  \end{equation}

  and
  \begin{equation}
    \{***\} \equiv \{ BX_1B^T + \tilde{A}X_3B^T + BX_2\tilde{A}^T + \tilde{A}X_4\tilde{A}^T \} \, .
  \end{equation}

  Expanding further considering the orders of $d$ contained in $[left]$, we can pull off the gradient contribution by looking at the first order terms in $d$:

  \begin{equation}
    \mathscr{D} f = tr(\{***\} M \{**\} + \{**\}M \{***\} )
  \end{equation}

  Now note that by making liberal use of the symmetry properties of the blocks of $X$, $\{***\}^T = \{***\}$. Similarly, $\{*\}^T = \{*\}$ therefore $\{**\}^T = \{**\}$. Recalling that traces are invariant under transpose, and that the trace of a sum is the sum of the traces, along with the fact that $M$ is symmetric, we can reduce this to the following compact form:

  \begin{equation}
    \mathscr{D} f = 2 tr(\{***\} M \{**\} )
  \end{equation}

  A remarkable simplification is that $\{***\} = \rv{B}{\tilde{A}}X\rv{B}{\tilde{A}}^T = \pi_V$. Expanding out the full expression and utilizing the cyclic property of the trace to push all of the $d$ and $d^T$ terms to the right, we find the following expression:
  \begin{eqnarray}
    \mathscr{D} f &=& 2 tr( X_1B^T\pi_V M d + \pi_VMBX_1d^T+\pi_VM\tilde{A}X_3d^T+X_2\tilde{A}^T\pi_vMd-X_1B^T\pi_VMBX_1B^Td \nonumber\\
    &-& BX_1B^T\pi_VMBX_1d^T - \tilde{A}X_3B^T\pi_VMBX_1d^T-X_1B^T\pi_VM\tilde{A}X_3B^Td -BX_1B^T\pi_VM\tilde{A}X_3d^T \nonumber\\
    &-& \tilde{A}X_3B^T\pi_VM\tilde{A}X_3d^T - X_1B^T\pi_VMBX_2\tilde{A}^Td - X_1B^T\pi_VM\tilde{A}X_4\tilde{A}^Td - X_2\tilde{A}^T\pi_VMBX_1B^Td \\
    &-& BX_2\tilde{A}^T\pi_VMBX_1d^T - \tilde{A}X_4 \tilde{A}^T \pi_VMBX_1d^T - X_2\tilde{A}^T \pi_VM\tilde{A}X_3B^Td - BX_2\tilde{A}^T\pi_VM\tilde{A}X_3d^T \nonumber\\
    &-& \tilde{A}X_4\tilde{A}^T\pi_VM\tilde{A}X_3d^T - X_2\tilde{A}^T\pi_VMBX_2\tilde{A}^Td - X_2\tilde{A}^T\pi_VM\tilde{A}X_4 \tilde{A}^T d)\nonumber
  \end{eqnarray}

  Using the trace to relate to the norm on our space, we can rewrite these terms to remove all transposes on the perturbation $d$ and read off the effect of the gradient:

  \begin{eqnarray}
    \mathscr{D} f &=& 2 ( \langle M \pi_V B X_1 ,d\rangle + \langle \pi_V M B X_1 ,d\rangle + \langle \pi_VM\tilde{A}X_3,d\rangle + \langle M\pi_V\tilde{A}X_3,d\rangle \nonumber \\
    &-& \langle BX_1B^TM\pi_VBX_1,d\rangle - \langle BX_1B^T\pi_VMBX_1,d\rangle - \langle \tilde{A}X_3B^T\pi_VMBX_1,d\rangle - \langle BX_2\tilde{A}^TM\pi_VBX_1,d\rangle \nonumber \\
    &-& \langle BX_1B^T\pi_VM\tilde{A}X_3,d\rangle - \langle \tilde{A}X_3B^T\pi_VM\tilde{A}X_3,d\rangle - \langle \tilde{A}X_3B^TM\pi_VBX_1 ,d\rangle - \langle \tilde{A}X_4\tilde{A}^T M\pi_VBX_1 ,d\rangle \nonumber \\
    &-& \langle BX_1B^TM\pi_V\tilde{A}X_3,d\rangle - \langle BX_2\tilde{A}^T\pi_VMBX_1,d\rangle - \langle \tilde{A}X_4\tilde{A}^T\pi_VMBX_1,d\rangle - \langle BX_2\tilde{A}^TM\pi_V\tilde{A}X_3,d\rangle  \\
    &-& \langle BX_2\tilde{A}^T\pi_VM\tilde{A}X_3,d\rangle - \langle \tilde{A}X_4\tilde{A}^T\pi_VM\tilde{A}X_3,d\rangle - \langle \tilde{A}X_3B^TM\pi_V\tilde{A}X_3,d\rangle - \langle\tilde{A}X_4\tilde{A}^TM\pi_V\tilde{A}X_3 ,d\rangle \nonumber )
  \end{eqnarray}

  With a clever observation, one can rearrange these terms back into a matrix equation for the gradient itself. Using $\mathscr{S}(\cdot)$ to denote the symmetric part of a square matrix given by $\mathscr{S} \, : \, M \mapsto \frac{M + M^T}{2}$ , we finally produce the following equation:

  \begin{equation}
    \nabla_B f = ( \mathbb{I}_{n\times n} - \pi_V) \mathscr{S}(M\pi_V) A X \cv{\mathbb{I}_{k\times k}}{\mathbb{O}_{i\times k}}
    \end{equation}

    A strong intuition might suggest that the gradient with respect to $\tilde{A}$ has a similar form with the other projected column of $X$. This can be checked to be the case:

    \begin{equation}
      \nabla_{\tilde{A}} f = ( \mathbb{I}_{n\times n} - \pi_V) \mathscr{S}(M\pi_V) A X \cv{\mathbb{O}_{k\times i}}{\mathbb{I}_{i\times i}}
    \end{equation}

\subsubsection{$\nabla_{\cdot} g$}

 Our computation for the second term $g(B,N)$ is much simpler.

 We can immediately see that

 \begin{equation}
   \nabla_{\tilde{A}} g = 0 \, .
 \end{equation}

 For $\nabla_B g$, we take a similar approach to the above computation, beginning with our expansion in $d$ up to linear order. For shorthand, we let $Y \equiv (B^TB)^{-1}$.

 \begin{eqnarray}
   g(B+d,N) &=& tr([(B+d)((B+d)^T(B+d))^{-1}(B+d)^T ] N [left]) \nonumber \\
    &=& tr([(B+d)(B^TB+B^Td+d^TB)^{-1}(B^T+d^T)] N [left]) + \mathcal{O}^+ \nonumber \\
    &=& tr([(B+d)(Y - Y(B^Td + d^TB)Y)(B^T+d^T)] N [left]) + \mathcal{O}^+  \\
    &=& tr([BYB^T + BY d^T - \pi_U dYB^T- BYd^T\pi_U + dYB^T] N [left]) + \mathcal{O}^+ \nonumber
   \end{eqnarray}

   Pulling off the first order term and exploiting trace and transpose symmetries as above, we find

   \begin{equation}
     \mathscr{D}g = 2tr( \{\diamond \} N \{\diamond \diamond \} ) \, ,
     \end{equation}

     where

   \begin{eqnarray}
    \{\diamond \} &\equiv&  BYB^T  \\
    \{\diamond \diamond \} &\equiv& BYd^T-\pi_UdYB^T - BYd^T\pi_U + dYB^T \, .
     \end{eqnarray}

  Massaging this into the form for our inner product $\langle \cdot, d \rangle$, we simplify to this expression:

  \begin{equation}
     \mathscr{D}g = 2 ( \langle NBY,d \rangle - \langle \pi_UNBY,d \rangle)
  \end{equation}

  We can finally read off our last gradient as follows:

  \begin{equation}
     \nabla_B g = 2 (\mathbb{I}_{n\times n} - \pi_U) N B Y
  \end{equation}

  \subsection{Hessians}

  At this stage, let the perturbed coordinates in a secondary direction be distinguished by a dot notation, $\dot X$.

  Handling $g$ first this time, we note through the chain rule and derivatives of matrix inverses that

  \begin{eqnarray}
    \dot Y &=& Y (\dot B^T B + B^T \dot B)Y \nonumber \\
    \dot \pi_U &=& \dot B Y B^T + B \dot Y B^T +  B Y \dot B^T
  \end{eqnarray}

  This results in a Hessian of the form:

  \begin{equation}
    H(g) = 2 [(\mathbb{I}_{n\times n} - \pi_U) N ( \dot B Y + B \dot Y) -\dot \pi_U N B Y  ]
  \end{equation}

  Similarly, using
  \begin{eqnarray}
    \dot A &=& \rv{\dot B}{\dot{\tilde{A}}} \nonumber \\
    \dot X &=& X(\dot A^T A + A^T \dot A ) X \\
    \dot \pi_V &=& \dot A X A^T + A \dot X A^T + A X \dot A^T \, ,\nonumber
  \end{eqnarray}

  we find that our two components of the Hessian have the following form:

  \begin{equation}
    H(f) = [(\mathbb{I}_{n\times n} - \pi_V)(\mathscr{S}(M \dot\pi_V)AX + \mathscr{S}(M\pi_V)(\dot A X + A \dot X)) - \dot \pi_V  \mathscr{S}(M\pi_V) A X ] \{ \cv{\mathbb{O}_{k\times k}}{\mathbb{I}_{i\times k}} , \cv{\mathbb{I}_{k\times i}}{\mathbb{O}_{i\times i}} \} \, ,
  \end{equation}
  where the projection vector in curly braces is chosen according to which component of the gradient is being used.

\bibliography{MyLibrary.bib}

\end{document}